\documentclass{amsart}
\usepackage{amsmath}
\usepackage{amsfonts}
\usepackage{amssymb}
\usepackage{graphicx}
\usepackage{multirow}
\usepackage{stix}
\usepackage{courier}
\usepackage{wrapfig}
\usepackage{eucal}
\usepackage{url}
\numberwithin{equation}{section}
\begin{document}

\title[Quadri-Figures in Cayley-Klein Planes: All around the Newton Line]{Quadri-Figures in Cayley-Klein Planes:\\ All around the Newton Line}
\author{Manfred Evers}
\curraddr[Manfred Evers]{Bendenkamp 21, 40880 Ratingen, Germany}
\email[Manfred Evers]{manfred\_evers@yahoo.com}
\date{\today}

\begin{abstract}
The Newton line and the associated theorems by Newton and Gauss for tetragons and quadrilaterals are closely linked to some other theorems of Euclidean geometry: a theorem by B\^ocher on the existence of a nine-point conic of a quadrangle, a theorem by Shatunov and Tokarev, and a theorem by Anne. This paper examines to which extent all these theorems can be transferred to other metric planes, in particular the elliptic and hyperbolic planes. \end{abstract}

\maketitle \hspace*{\fill}\vspace*{-1 mm}\\

\section*{The Newton line of a tetragon}
Consider a tetragon $\mathbf{T} = ABCD$ in the Euclidean plane. 
If this tetragon is a parallelogram, then the centers $M_{AC}, M_{BD}$ of the diagonal segments $[A,C]$ and $[B,D]$ both agree with the centroid of the quadrangle $\mathbf{QA} =\{A, B, C, D\}$, a point which we can get as the intersection of the two lines that connect the midpoints of opposite sides of the tetragon. If the tetragon is not a parallelogram, then the centroid and the points $M_{AC}, M_{BD}$ are collinear on a line which is called the  \textit{Newton line} of the tetragon. Newton discovered\footnote{$^)$ See \cite{HS} for historical background.}$^)$ that if there is a circle which touches the four sidelines of $\mathbf{T}$ (in this case $\mathbf{T}$ is usually called a
\textit{tangential quadrilateral}), then $\mathbf{T}$ is either a rhombus or the center of this circle is a point on a line which is now called \textit{Newton line}. 

\begin{figure}[!htbp]
\includegraphics[height=6cm]{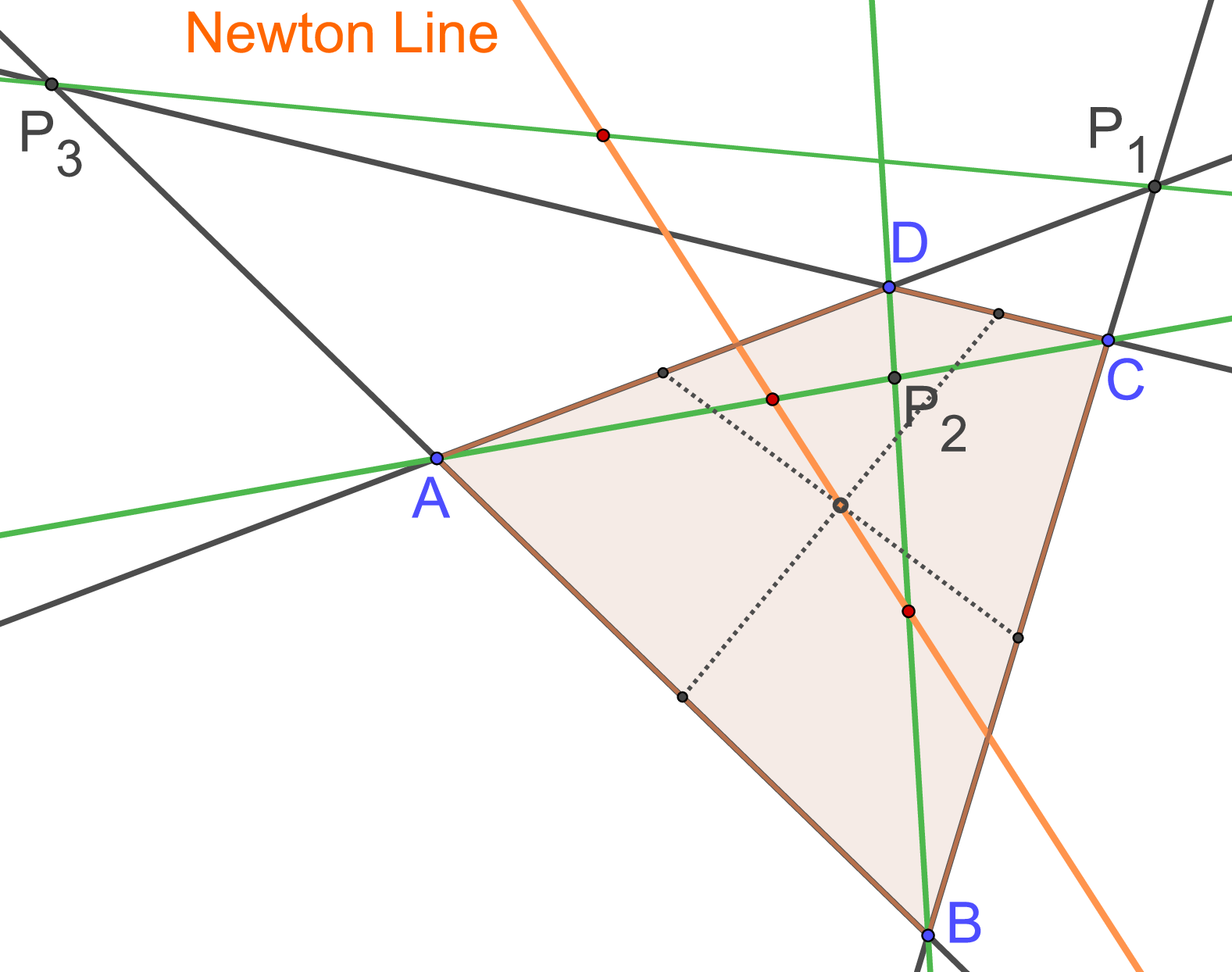}
\caption{All figures were created with the software program \textit{GeoGebra}.}
\vspace*{1 mm}\end{figure}

We expand the Euclidean plane to a metric projective plane by adding a line at infinity $\mathcal{L}^\infty$. 
Each tetragon $\mathbf{T}$ can be uniquely assigned a complete quadrangle $\mathbf{QA}$ and a complete quadrilateral $\mathbf{QL}$. The complete quadrangle $\mathbf{QA}$ has the same vertices as $\mathbf{T}$, but six sidelines; in addition to the sides of $\mathbf{T}$, the lines $A{\,\vee\,}C$ and $B{\,\vee\,}D$ are also sides of $\mathbf{QA}$. The complete quadrilateral $\mathbf{QL}$ has the same sidelines as $\mathbf{T}$, but six vertices; in addition to the points $A, B, C, D$, the intersections of opposite sidelines of $\mathbf{T}$ are also vertices of  $\mathbf{QL}$.\\
The points $P_1 = (A{\,\vee\,}D)\wedge(B{\,\vee\,}C) , P_3=(C{\,\vee\,}D)\wedge(A{\,\vee\,}B)$ together with the point $P_2=(B\vee D)\wedge(A\vee C)$ are the \textit{diagonal points} of $\mathbf{QA}$, the lines $P_1{\,\vee\,}P_3, A{\,\vee\,}C, B{\,\vee\,}D$ are the \textit{diagonals} of $\mathbf{QL}$. The diagonal points and the diagonals can be directly assigned to the tetragon $\mathbf{T}$.\\ 
%We put $\overline{\mathbf{QA}}{\,{{\scriptstyle{:}}\!\!=}\,}\{A,B,C,D,P_1,P_2,P_3\}$, $\partial\,\overline{\mathbf{QA}}{\,{{\scriptstyle{:}}\!\!=}\,}\overline{\mathbf{QA}}{\,\smallsetminus}{\hspace{-4.3pt}\smallsetminus\,}\mathbf{QA} = \{P_1,P_2,P_3\}$.\\

Gauss discovered that the Newton line of $\mathbf{T}$ goes not only through the midpoints of the segments $[A,C]$ and $[B,D]$, but also through the middle of the segment $[P_1,P_3]$. For this reason, this straight line is often referred to as the \textit{Gauss-Newton line} of the complete quadrilateral $\mathbf{QL}$ and is uniquely determined. 
Gauss' discovery does not only apply to the Euclidean plane, but - as can be easily checked - to every metric affine plane, which we always imagine to be embedded in a (metric) projective plane. 

In this projective plane we can assign each line-segment not only one midpoint but two: an inner midpoint, which lies within this segment, and an outer midpoint, which lies on the line $\mathcal{L}^\infty$ outside the segment.

Thus, Gauss' observation (Gauss' theorem) can also be formulated as follows: 
The six midpoints of the diagonals of a complete quadrilateral $\mathbf{QL}$ lie on a singular conic.\\

One is tempted to check to what extent this statement applies to other Cayley-Klein geo\-metries.   Experiments with mathematical geometry software (\textit{Geogebra} \cite{GG}) suggest that the six midpoints of the diagonals of a complete quadrilateral in an elliptic plane always lie on a single conic and that the six semi-midpoints of a complete quadrilateral in a hyperbolic plane\footnote{$^)$ When we use the term "hyperbolic plane", we always mean the extended hyperbolic plane.}$^)$
also lie on a single conic.\\
These observations are confirmed by \vspace*{-0.6 mm}\\

\begin{figure}[!htbp]
\includegraphics[height=8cm]{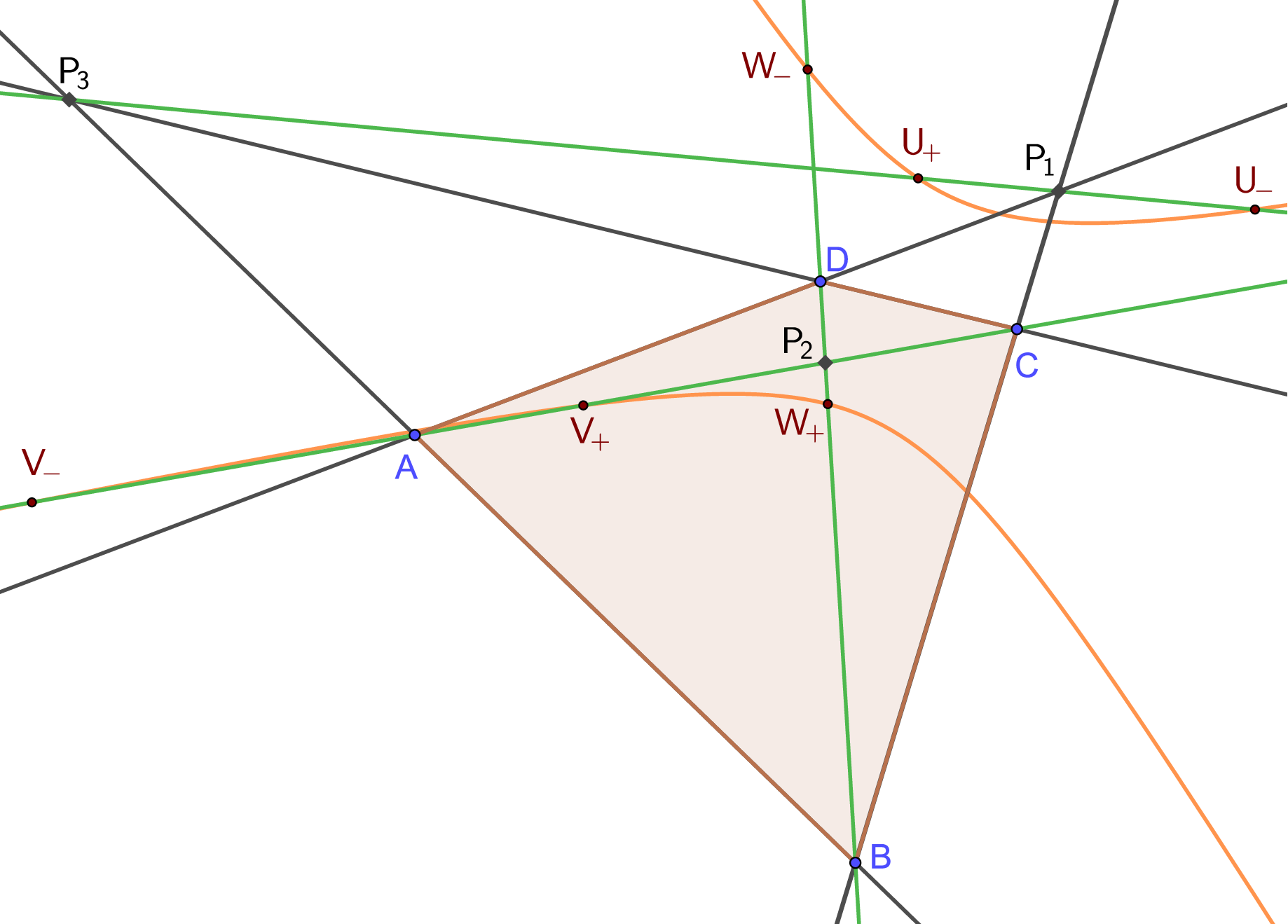}
\caption{}
\vspace*{1 mm}\end{figure}

\noindent\textit{Theorem} 1.
Let $A, B, C, D$ be the vertices of a tetragon in the projective plane $\textrm{P}(\mathbb{R}^3)$, let $P_1, P_2, P_3$ be its diagonal points and $A\vee C, B\vee D, P_1\vee P_3$ its diagonals.    \\
Whenever we choose six distinct points $U_+, U_-, V_+, V_-, W_+, W_-$ such that the three quadruples $(U_-, P_1, U_+,P_3)$, $(V_-,A,V_+,C)$, $(W_-,B,W_+,D)$ form harmonic ranges, 
then there exists a conic passing through all these six points.\vspace*{0.7 mm}\\
\textit{Proof}.  The points  $P_1,P_2,P_3$ are not collinear and form a generating
system for the projective plane. If a point $Q$ has homogenous barycentric coordinates $q_1, q_2, q_3$with respect to the point triple $(P_1,P_2,P_3)$, we write $Q = [q_1{:}q_2{:}q_3]$. Thus, $P_1=[1{:}0{:}0], 
P_2=[0{:}1{:}0], P_3=[0{:}0{:}1]$. For the points $A, B, C, D$ there exist nonzero real numbers $s,t$ such that $A=[-s{:}1{:}t]$, $B=[-s{:}1{:}{-t}]$, $C=[s{:}1{:}{-t}]$, $D=[s{:}1{:}t]$, and we can find  real numbers $u,v,w$ with $0<u,v,w<1$ and $U_\pm=[1{\,:\,}0{\,:\,}\pm u] , V_{\pm}=[(-1\pm v)s{\,:\,}1\pm v{\,:\,}(1\mp v)t],W_{\pm} =[(-1\pm w)s{\,:\,}1\pm w{\,:\,}(-1\pm w)t]$.
It can be easily varified that $U_\pm,V_\pm,W_\pm$ are points on a conic
with matrix \vspace*{1.5 mm}\\
\centerline{$\displaystyle\mathfrak{M}{\,=\,} \left(\begin{array}{ccc} 
\!s t u^2 (v^2{-}1) (1{-}w^2) &\!\! t (s^2 u^2{-}t^2) (v^2 w^2{-}1) &\!\! 0\\ 
\!t (s^2 u^2{-}t^2) (v^2 w^2{-}1) &\!\! s t (v^2{-}1) (1{-}w^2) (s^2 u^2{-}t^2) &\!\! s (v^2{-}w^2) (s^2 u^2{-}t^2)\\ 
\!0 &\!\!s (v^2{-}w^2) (s^2 u^2{-}t^2)&\!\! s t (1{-}v^2) (1{-}w^2)
\end{array}\right).\,\;\Box$}\vspace*{1.2 mm}\\

According to the principle of duality in projective geometry, the dual version of the Theorem 1 is also a valid theorem.\vspace*{1.7 mm}\\
\textit{Theorem} 2. Let $A, B, C, D$ be the vertices of a tetragon in a Euclidean plane $\textrm{P}(\mathbb{R}^3)$ and let $P_1, P_2, P_3$ be the diagonal points.\\
Whenever we choose six distinct lines $\mathcal{U}_+,\mathcal{U}_-,\mathcal{V}_+,\mathcal{V}_-,\mathcal{W}_+,\mathcal{W}_-$ such that the three quadruples $(A{\,\vee\,}D,\mathcal{U}_-,B{\,\vee\,}C,\mathcal{U}_+)$, $(B{\,\vee\,}D,\mathcal{V}_-,A{\,\vee\,}C ,\mathcal{V}_+)$, $(C{\,\vee\,}D,\mathcal{W}_-,A{\,\vee\,}B,\mathcal{W}_+)$ form harmonic pencils, there exists a conic touching all these six lines $\mathcal{U}_+,\mathcal{U}_-,\mathcal{V}_+,\mathcal{V}_-,\mathcal{W}_+,\mathcal{W}_-$. (Definition: A line $\mathcal{L}$ \textit{touches} a conic $\mathcal{K}$ if $\mathcal{L}$ and $\mathcal{K}$ have a single common point.)\vspace*{2.2 mm}

As a consequence  we get the two following corollaries:\vspace*{1.1 mm}\\
\textit{Corollary} 1. 
Let $A, B, C, D$ be the vertices of a tetragon in a Euclidean plane, let $P_1, P_2, P_3$ be the diagonal points and $A\vee C, B\vee D, P_1\vee P_3$ the diagonals. Then there exists a conic touching all bisectors of the angles $\angle(BP_1 D)$, $\angle(CP_2D)$, $\angle(AP_3D)$  (these are altogether six lines). See Figure 3.\vspace*{-1 mm}\\
\begin{figure}[!hptb]
\begin{centering} 
\includegraphics[height=7cm]{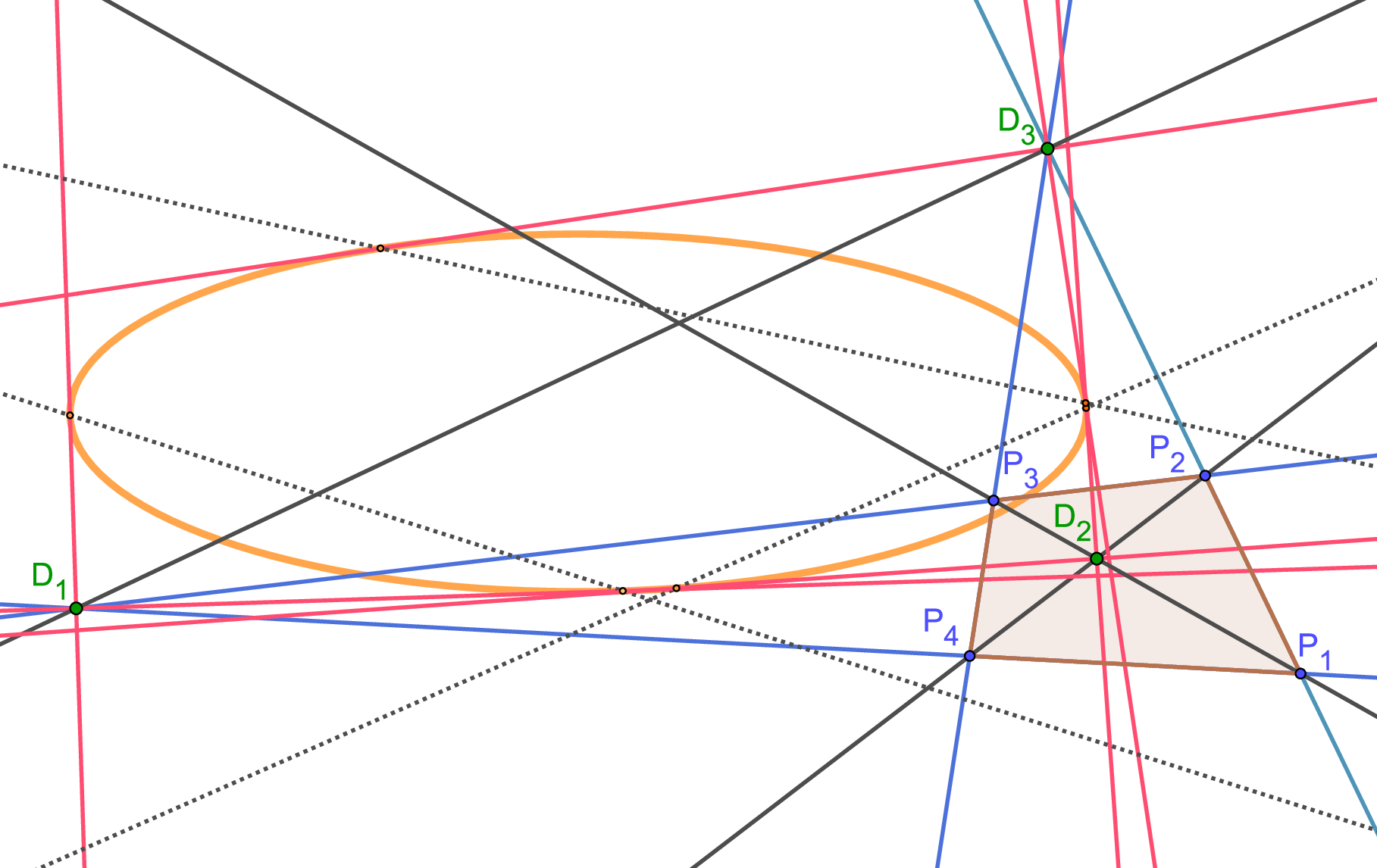}
\caption{}

\includegraphics[height=8cm]{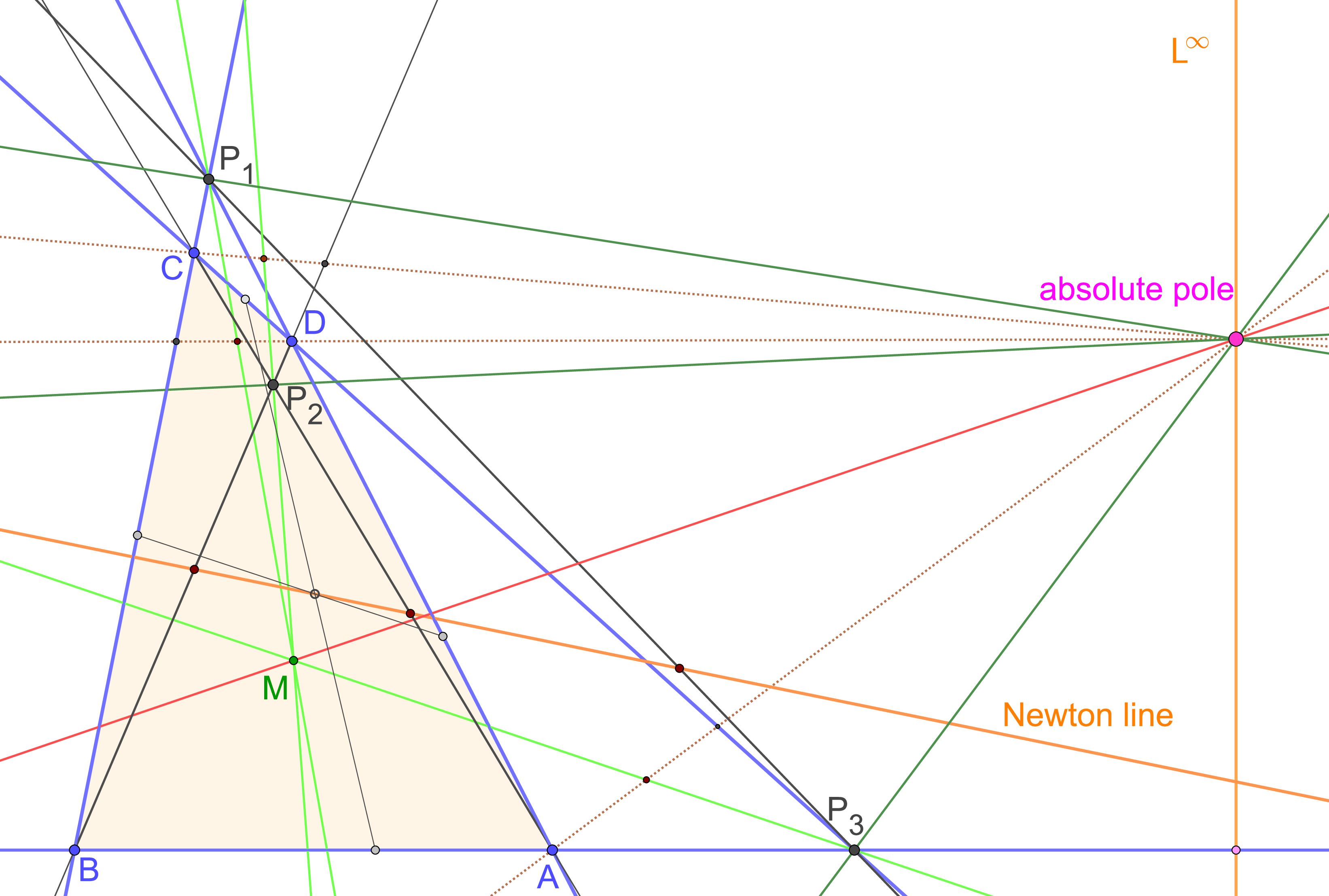}
\caption{}
\end{centering} 
\vspace*{0 mm}\end{figure}\vspace*{0.6mm}\\
\textit{Corollary} 2. 
Let $A, B, C, D$ be the vertices of a tetragon in a Galilean plane, let $P_1, P_2, P_3$ be the diagonal points and $A\vee C, B\vee D, P_1\vee P_3$ the diagonals. Then one of the two bisectors of each of the angles $\angle(AP_1 D)$, $\angle(AP_2C)$, $\angle(AP_3B)$ passes through the absolute pole $P_{abs}$, while the other angle bisectors meet at some point $M$ other than $P_{abs}$. See Figure 4.\vspace*{-3 mm}\\

\centerline{--------}\vspace*{-2mm}
\centerline{-----}\vspace*{0mm}
\centerline{\textit{Excursus}. Geometry on elliptic and hyperbolic planes.} \vspace*{3 mm}

We assume that the reader is familiar with projective, elliptic and hyperbolic geometry, but in order to introduce the terminology and fix notations, we give some basic definitions, rules and theorems. \vspace*{-1.4 mm}\\

%In the following we focus on elliptic and a hyperbolic projective planes. In order to formulate and prove theorems concerning tetragons, quadrugons and quadrilaterals in an elliptic and in a hyperbolic plane, we provide necessary foundations and terminology.\vspace*{-1.4 mm}\\

As already introduced, let points  $P_1=[1{:}0{:}0], P_2=[0{:}1{:}0], P_3=[0{:}0{:}1]$ generate the projective plane.\vspace*{-1.4 mm}\\ 

Poles and polars.
Every symmetric regular matrix $\mathfrak{M} = (\mathfrak{m}_{ij})_{1\le i,j\le 3} \in \mathbb{R}^3\times \mathbb{R}^3$ induces \\
\;-\;  a bilinear function $... {\scriptstyle{[\mathfrak{M}]}}...{:}\mathbb{R}^3\times \mathbb{R}^3\to \mathbb{R}$, $(p_1,p_2,p_3){\scriptstyle{[\mathfrak{M}]}}(q_1,q_2,q_3){\,{{\scriptstyle{:}}\!\!=}\,}\sum\limits_{1\le i,j\le 3} p_i \mathfrak{m}_{ij} q_j$\vspace*{-0.8 mm}\\
and\\
\;-\; a polarity between points and lines in a projective plane: For any point $P{\,{{\scriptstyle{:}}\!\!=}\,}[p_1{:}p_2{:}p_3]$ in a projective plane, the point set $P^{\mathfrak{M}}{\,{{\scriptstyle{:}}\!\!=}\,}\{[x_1{:}x_2{:}x_3]\,|\,$$(x_1,x_2,x_3){\scriptstyle{[\mathfrak{M}]}}(p_1,p_2,p_3)=0\}$ is a line which is called the \textit{polar} (\textit{line}) of $P$ with respect to $\mathfrak{M}$. For any line $\mathcal{L}$ in this projective plane, there exists a nonzero vector $(l_1,l_2,l_3)\in \mathbb{R}^3$ such that $\mathcal{L}=\{[x_1{:}x_2{:}x_3]\, |\, l_1 x_1+l_2 x_2+l_3 x_3 =0\}$ and there exists exactly one point $P=[p_1{:}p_2{:}p_3]$ in the plane satisfying $(p_1,p_2,p_3)=(l_1,l_2,l_3)\mathfrak{M}^{-1}$; this point is the \textsl{pole} of $\mathcal{L}$ with respect to $\mathfrak{M}$.\\

The metric structure.\\
The metric structure of the elliptic and the hyperbolic plane is given by a special regular symmetric real matrix $\mathfrak{G} =  (\mathfrak{g}_{ij})_{1\le i,j\le 3}$. The polarity defined by $\mathfrak{G}$ is called \textit{duality}. \\
For an elliptic plane, the matrix $\mathfrak{G}$ is definit and we may assume that it is positive definit and that $\mathfrak{g}_{11}=\mathfrak{g}_{22}=\mathfrak{g}_{33} = 1$. In the hyperbolic case $\mathfrak{G}$ is indefinite, so $\mathfrak{G}$ is the matrix of a  nonsingular conic, called the \textit{absolute conic}. Points on this conic are called \textit{isotropic}, all the others are \textit{anisotropic}. All points in the elliptic plane are anisotropic.\\
The elliptic and the hyperbolic plane are \textit{regular} Cayley-Klein planes, since the matrix $\mathfrak{G}$ is regular. The metric affine planes and their duals are \textit{singular} Cayley-Klein planes.\\

Isometries in regular Cayley-Klein planes.\\
Reflections: For any anisotropic point $M=[m_1{:}m_2{:}m_3]$ and any point $P=[p_1{:}p_2{:}p_3]$, 
the image of $P$ under a \textit{reflection} in $M$ is the point $Q = [q_1{:}q_2{:}q_3]$ with
\begin{equation*}
(q_1,q_2,q_3) = 2\frac{(m_1,m_2,m_3){\,\scriptscriptstyle{[\mathfrak{G}]}\,}(p_1,p_2,p_3)}{(m_1,m_2,m_3){\,\scriptscriptstyle{[\mathfrak{G}]}\,}(m_1,m_2,m_3)} (m_1,m_2,m_3) - (p_1,p_2,p_3).
\end{equation*} \\
Reflections are involutions, mapping isotropic points to isotropic points and collinear points to collinear points. The reflection in $M$ leaves this point $M$ fixed, but also all points on the dual line of $M$.\\
A composition of finitely many reflections is an \textit{isometry}.\\
If $\mathcal{S}_1, \mathcal{S}_2$ are sets in an elliptic or in a hyperbolic plane, $\mathcal{S}_1$ is 
\textit{congruent to} $\mathcal{S}_2$ (we write $\mathcal{S}_1 \cong \mathcal{S}_2$) if there exists an isometry which maps  $\mathcal{S}_1$ onto $\mathcal{S}_2$.
Obviously, congruence $\cong$ is an equivalence relation.\vspace*{1 mm}\\
While in an elliptic plane $\{P\} \cong \{Q\}$ holds for all points $P$ and $Q$, this is not the case in a hyperbolic plane; here $\{P\} \cong \{Q\}$ precisely when $P$ and $Q$ are both isotropic or both points lie inside or both outside the absolute conic.\vspace*{-1 mm} \\

Symmetry points of a set.\\
Let $\mathcal{S}$ is be subset of the elliptic or the hyperbolic plane. Then an anisotropic point $P$ is a
\textit{symmetry point} of $\mathcal{S}$, if the image of $\mathcal{S}$ under a reflection in $P$ is a subset of $\mathcal{S}$. A line $\mathcal{L}$ is a \textit{symmetry axis} of $\mathcal{S}$ if its dual point is a symmetry point of $\mathcal{S}$.\vspace*{1 mm} \\ 
We give two examples:\vspace*{1 mm} \\
(1) If $\mathcal{S}$ consists of a single point $P$, then all anisotropic points on its dual line are the symmetry points. If $P$ is anisotropic, then $P$ is also a symmetry point. 
(2) If $\mathcal{S}$ is a line, then each anisotropic point on $\mathcal{S}$ is a symmetry point. If $\mathcal{S}$ is an anisotropic line, then the pole of $\mathcal{L}$ is also a symmetry point.\\
%(3) If $\mathcal{S}$ is the union of two distinct lines which meet at an anisotropic point $P$, then the point $P$ and the duals of the two angle bisectors are the symmetry points.\\

Semi-midpoints and line segments.\\
Define a function $\chi{:}\;V \rightarrow \{-1,0,1\}$ by\vspace*{-1 mm}  
\[\chi(p_1,p_2,p_3) = 
\begin{cases}
 \;\;\,0,\;\text{if}\,(p_1,p_2,p_3) = (0,0,0)\;,\\
 \;\;\,1,\;\text{if}\,(p_1,p_2,p_3) > (0,0,0)\; \textrm{with respect to the reverse lexicographic order,} \\
-1,\;\text{if}\, (p_1,p_2,p_3) < (0,0,0)\; \textrm{with respect to the reverse lexicographic order.} \\
\end{cases}\vspace*{-0.9 mm}\]
For an anisotropic point $P = [p_1{:}p_2{:}p_3]$ let $P^\circ$ be the vector \vspace*{-2.3 mm}\\
\[\,P^\circ\!{\,{{\scriptstyle{:}}\!\!=}\,}\frac{\chi(p_1,p_2,p_3)}{\sqrt{|(p_1,p_2,p_3){\,\scriptscriptstyle{[\mathfrak{G}]}\,}(p_1,p_2,p_3)|}} (p_1,p_2,p_3) \in \mathbb{R}^3.\; \\
\] \vspace*{0.8mm}
\noindent 
%Given two distinct anisotropic points $P$ and $Q$, we call the point $S=[s_1{:}s_2{:}s_3]$ with $(s_1,s_2,s_3)=P^\circ + Q^\circ$ the \textit{inner semi-midpoint} of $P$ and $Q$ and the point $T=[t_1{:}t_2{:}t_3]$ with $(t_1,t_2,t_3)=P^\circ - Q^\circ$ the \textit{outer semi-midpoint} of $P$ and $Q$.\\ The four points $P, S, Q, T$ form a harmonic range.\\
Given two distinct anisotropic points $P$ and $Q$, we call the point $M=[m_1{:}m_2{:}m_3]$ with $(m_1,m_2,m_3)=P^\circ + Q^\circ$ and the point $N=[n_1{:}n_2{:}n_3]$ with $(n_1,n_2,n_3)=P^\circ - Q^\circ$ semi-midpoints of $P$ and $Q$. The four points $P, M, Q, N$ form a harmonic range.
If the mirror image of $P$ in $M$ is $Q$, then $M$ and $N$ are called (\textit{proper}) \textit{midpoints}, otherwise they are called \textit{pseudo-midpoints}. In an elliptic plane all semi-midpoints are proper midpoints. While proper midpoints are always anisotropic, pseudo-midpoints do not have to be.\\
For any two distinct points $P$ and $Q$ we introduce two line segments $[P,Q]_{+}, [P,Q]_{-}.\;$ These are closures of the two connected components of the set $P{\vee}Q \;{\smallsetminus}{\hspace{-4.3pt}\smallsetminus}\;\{P, Q\}$. $[P,Q]_{+}$ is characterized by the property that whenever $R$ and $S$ are distinct anisotropic points in $[P,Q]_{+}$, then $M=[m_1{:}m_2{:}m_3]$ with $(m_1{,}m_2{,}m_3) = R^\circ+S^\circ$ is also a point in $[P,Q]_{+}$. If $R$ and $S$ are distinct anisotropic points in $[P,Q]_{-}$, then $N=[n_1{:}n_2{:}n_3]$ with $(n_1{,}n_2{,}n_3) = R^\circ-S^\circ$  also belongs to $[P,Q]_{-}$.\\ 
We call a semi-midpoint of a line segment \textit{inner semi-midpoint} if it lies within the segment, otherwise \textit{outer semi-midpoint}.

Perpendicular bisector of a line segment.\\
Suppose $[A,B]_\pm$ is a segment which has a proper inner midpoint $M_\pm$. Then there is the line which passes through this midpoint and through the dual of $A\vee B$; it is called the \textit{perpendicular bisector} of $[A,B]_\pm$. This bisector is the dual of the point $M_\mp$.\\

Triangles.\\
Given three non-collinear anisotropic points $A,B,C$. Then the closure of each connected component of the set $\textrm{P}(\mathbb{R}^3){\,\smallsetminus}{\hspace{-4.3pt}\smallsetminus\,}((A\vee B)\cup (B\vee C) \cup (A\vee A))$ is a \textit{triangle} with vertices $A,B,C$. There are four such triangles $\Delta_i(A,B,C), i = 0,1,2,3$; their bounderies are\\
\centerline{$\partial \Delta_0 = [A,B]_+\cup[B,C]_+\cup [C,A]_+$\;,\;$\partial \Delta_1 = [A,B]_-\cup[B,C]_+\cup [C,A]_-$\,,}\\
\centerline{$\partial \Delta_2 = [A,B]_-\cup[B,C]_-\cup [C,A]_+$\;,\;$\partial \Delta_3 = [A,B]_+\cup[B,C]_-\cup [C,A]_-$\,.}\\

Conics.\\
Let $\mathfrak{M}$ be a symmetric matrix. If $\mathfrak{M}$ is indefinit, then there is a real conic $\mathcal{K}=\mathcal{K}(\mathfrak{M})$ which consists of all points $[p_1{:}p_2{:}p_3]$ with $(p_1,p_2,p_3){\scriptstyle{[\mathfrak{M}}]}(p_1,p_2,p_3)=0$. $\mathcal{K}$ is a nonsingular conic iff  $\mathfrak{M}$ is a nonsingular matrix. \\
Let $\mathcal{K}=\mathcal{K}(\mathfrak{M})$ be a nonsingular conic, but not the absolute.\\
Then an anisotropic point $P$ is a symmetry point of $\mathcal{K}$ iff $P^{\mathfrak{M}} = P^{\mathfrak{G}}$. \\
We call the conic $\{ P\;|\; P^{\mathfrak{G}}\;\mathrm{is\;a\;tangent\;line\;of\;}\mathcal{K}\}$  the \textit{dual} of $\mathcal{K}$. (According to this definition, the dual of a point conic is a point conic and not a line conic.) \\
For a complete classification of conics in elliptic and hyperbolic planes, see \cite{Iz, St}.\\ 

Circles.\\
Circles are special nonsingular real conics. (Double points and double lines are not considered circles here.)\\
If $\mathcal{K}$ is a nonsingular real conic (different from the absolute), then $\mathcal{K}$ is called a \textit{circle} if there is a line $\mathcal{L}$ whose anisotropic points are all symmetry points of $\mathcal{K}$. (If $\mathcal{K}$ is not a circle, $\mathcal{K}$ can have at most three symmetry points.)\\
The pole $M$ of $\mathcal{L}$ is called the center of the circle $\mathcal{K}$. Suppose $M$ is anisotropic, then $M$ is also a symmetry point and for any points  $P$ and $Q$ on the circle  the two segments $[M,P]_+$ and $[M,Q]_+$ are congruent.

If $A,B,C$ are three noncollinear anisotropic points, then there is a circle passing through these points iff
all segments $[A,B]_+, [B,C]_+, [C,A]_+$ have proper midpoints. But in this case there is not only one but there are four circles running through $A,B,C$. There is a 1 to 1 relation between these circles and the triangles $\Delta_i(A,B,C),\; i=0,1,2,3$. The three perpendicular bisectors of each of these triangles meet at a point, which is the center of exactly one of these circles. 

\begin{figure}[!htbp]
\includegraphics[height=8cm]{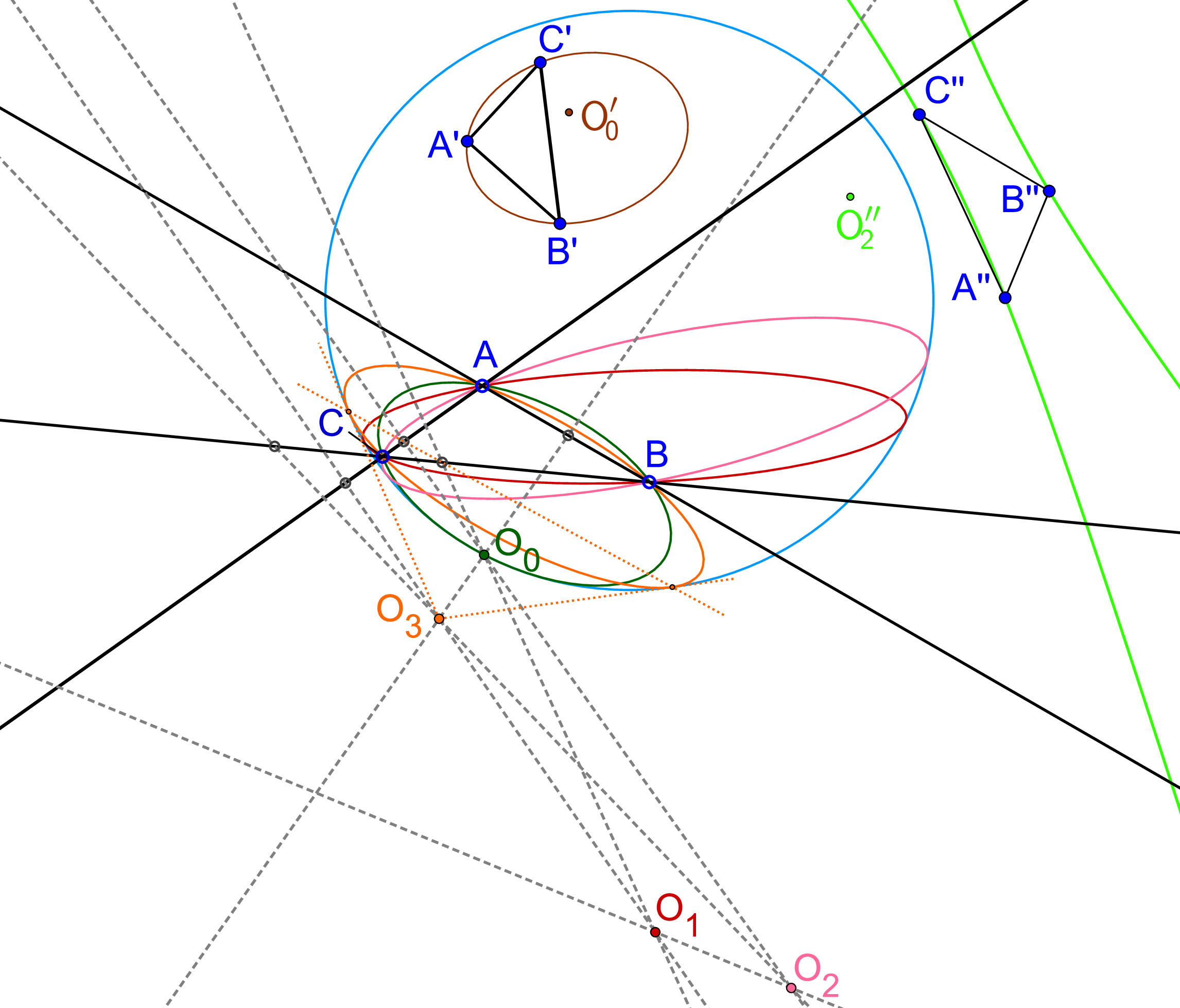}
\caption{Circles in the hyperbolic plane. The blue conic is the absolute. The four circumcircles of triangle $ABC$ are cycles, the one with center $O_0$ is a horocycle.
The brown and the light green circle are proper circles.}
\vspace*{1 mm}\end{figure}

While every circle in the elliptic plane can be considered as a special ellipse, the situation in the hyperbolic plane is more complicated. Figure 5 shows different types of circumcircles of triangles in the hyperbolic plane. A \textit{proper circle} is one whose center lies inside the circle. However, there are two types within this group: there are proper circles that lie inside the absolute conic and those that lie outside it.\\
Circles whose centers are not inside the circle are called \textit{cycles}. If the center lies outside the circle, the circle is called  \textit{hypercycle}; and if the center is a point on the circle, it is called \textit{horocycle}.
There are hypercycles and horocycles whose anisotropic points lie inside the absolute conic and those whose anisotropic points lie outside it.
While the centers of proper circles and hypercycles are symmetry points, the center of a horocycle is not, because it is an isotropic point. \vspace*{-3.7 mm}\\

\centerline{-----}\vspace*{-2mm}
\centerline{--------}\vspace*{2mm}

We end Section 1 with\vspace*{-2 mm}\\

\noindent\textit{Theorem} 3.
Let $A, B, C, D$ be the vertices of a tetragon in an elliptic or in a hyperbolic plane. Let $P_1, P_2, P_3$ be the diagonal points and $A\vee C, B\vee D, P_1\vee P_2$ the diagonals. Here and in the following we assume that the vertices and diagonal points of a tetragon are anisotropic points. \\ 
There exists a conic passing through the semi-midpoints of $P_1,P_3$, the semi-midpoints of $A,C$ and the semi-midpoints of $B,D$.\\
\textit{Proof}. Put $u{\,{{\scriptstyle{:}}\!\!=}\,}1$, $v{\,{{\scriptstyle{:}}\!\!=}\,}\sqrt{|C {\scriptstyle{[\mathfrak{G}]}} C|}/\sqrt{|A {\scriptstyle{[\mathfrak{G}]}} A|}$, $w{\,{{\scriptstyle{:}}\!\!=}\,}\sqrt{|D {\scriptstyle{[\mathfrak{G}]}} D|}/\sqrt{|B {\scriptstyle{[\mathfrak{G}]}} B|}$. Then the matrix $\mathfrak{M}$ as described in the proof of Theorem 1 is a matrix of a conic with the required properties. $\Box$

\section*{The nine-point conic of a quadrangle}
We adopt the definition of the nine-point conic in a metric affine plane as it is given by Chris van Tienhoven in \cite{Ti}. 
The \textit{nine-point conic} of a quadrangle in a metric-affine plane is the conic through the midpoints of all possible line segments connecting any two vertices of the quadrangle.\\ 
This definition initially only justifies the name \textit{six-point conic}. But apart from these six midpoints it also passes through the intersection points of all possible pairs of lines connecting the vertices of the quadrangle; these are another three points. \\
In the year 1892 Maxime B\^ocher described this conic in a paper \cite{Bo}.\footnote{$^)$The existence of such a conic was probably already known earlier, cf. \cite{Od}.}$^)$ \\

\noindent Projective version of the Nine-Point Conic \textit{Theorem} (cf. \cite{Al, Od}). Let $\mathbf{QA} =\{A, B, C, D\}$ be a quadrangle and let $\mathcal{L}$ be a line which meets the six side lines of the complete quadrangle in points $P_{12} \in A\vee B, P_{13}\in A\vee C,\dots,P_{34}\in C\vee D$. If $P'_{12}$ is the harmonic conjugate of $P_{12}$ with respect to $A$ and $B$, $P'_{13}$ the harmonic conjugate $P_{13}$ with respect to $A$ and $C$, ..., then there exists a conic passing through these six points $P'_{12},\dots,P'_{34}$ and through the diagonal points of $\mathbf{QA}$. We will call this conic the \textit{nine-point conic of} $(\mathbf{QA},\mathcal{L})$.\\

\noindent\textit{Theorem} (Odehnal \cite{Od}). We adopt the assumptions of the previous theorem.
In addition, let $\mathcal{K}$ be any conic passing through the points $A,B,C,D$ and let $P$ be the pole of $\mathcal{L}$ with respect to $\mathcal{K}$. 
Choose any two distinct points $Q_1, Q_2\in \{A,B,C,D\}$, then the line through $P$ and the pole of $Q_1\vee Q_2$ with respect to $\mathcal{K}$ meets the line $Q_1\vee Q_2$ in a point on the nine-point conic of $(\mathbf{QA},\mathcal{L})$. The point $P$ is also a point on this conic. See Figure 6.\\

\begin{figure}[!hbpt]
\includegraphics[height=6cm]{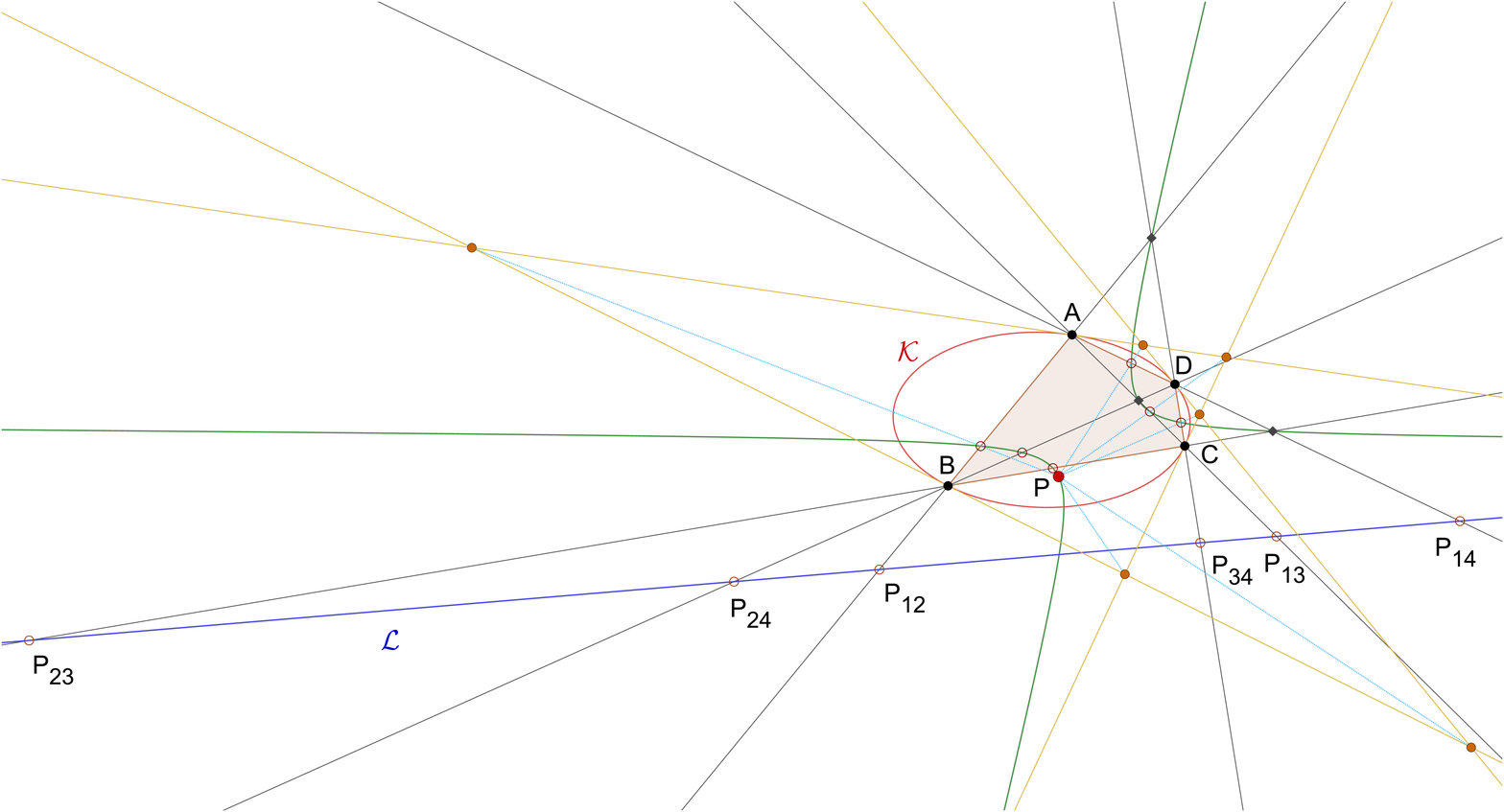}
\caption{}
\vspace*{1 mm}\end{figure} 
\begin{figure}[!hbpt]
\includegraphics[height=6cm]{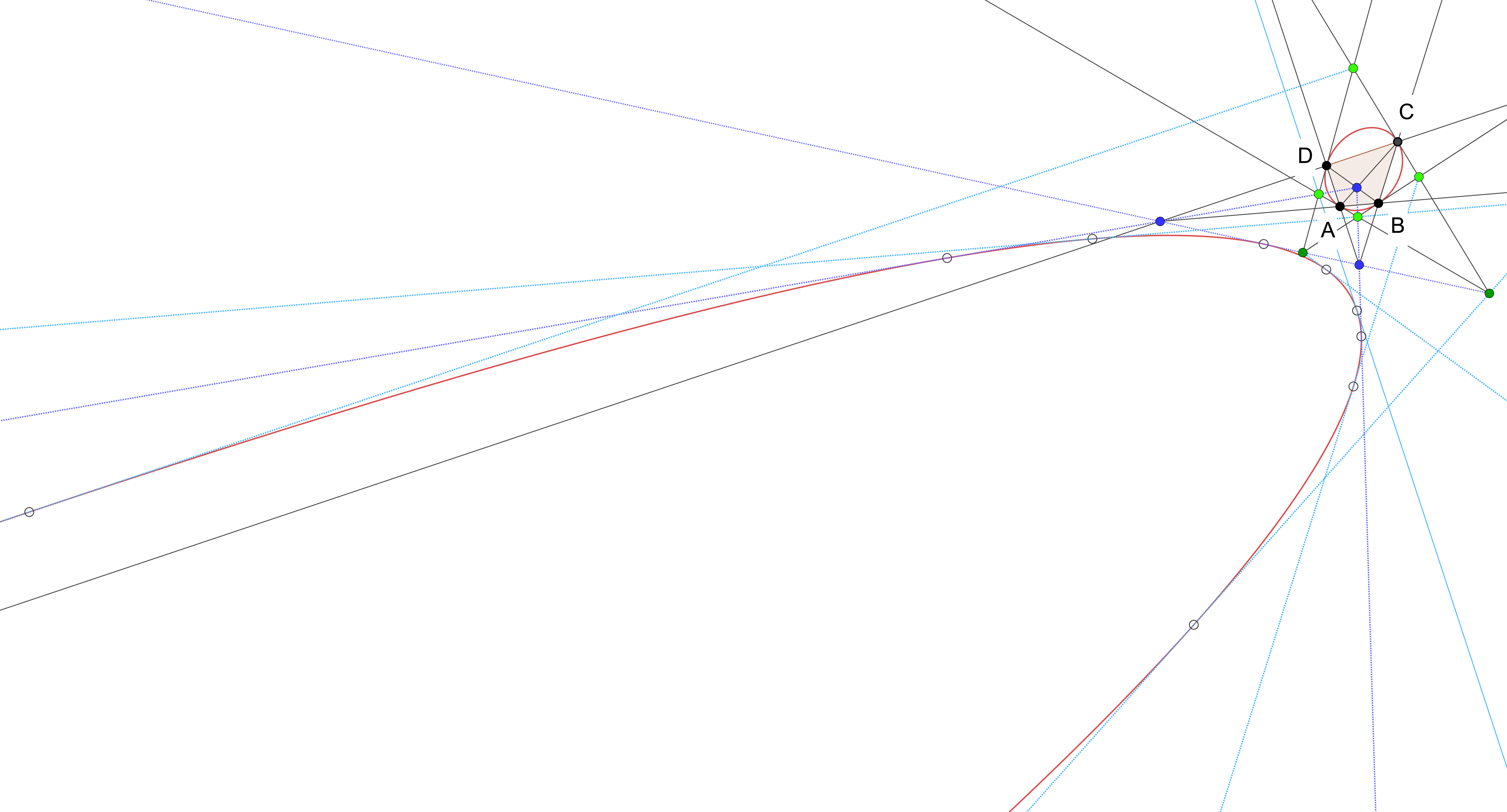}
\caption{}
\vspace*{1 mm}\end{figure}

\noindent\textit{Corollary} (A special dual version of the previous theorem). Let $\mathbf{QA} =\{A, B, C, D\}$ be a quadrangle in a metric affine plane and $\mathcal{K}$ be a conic passing through the points $A,B,C,D$. Then there is a parabola $\mathcal{P}$ with the following property:\\
Whenever we choose a line $\mathcal{L}$ which connects any two of the points $A,B,C,D$, the parallel of $\mathcal{L}$ through the pole of $\mathcal{L}$ wrt $\mathcal{K}$ is a tangent to $\mathcal{P}$. 
Furthermore, the lines that connect two diagonal points of $\mathbf{QA}$ are also tangents to $\mathcal{P}$. See Figure 7	\\

\begin{figure}[!htbp]
\begin{centering}
\includegraphics[height=7.6cm]{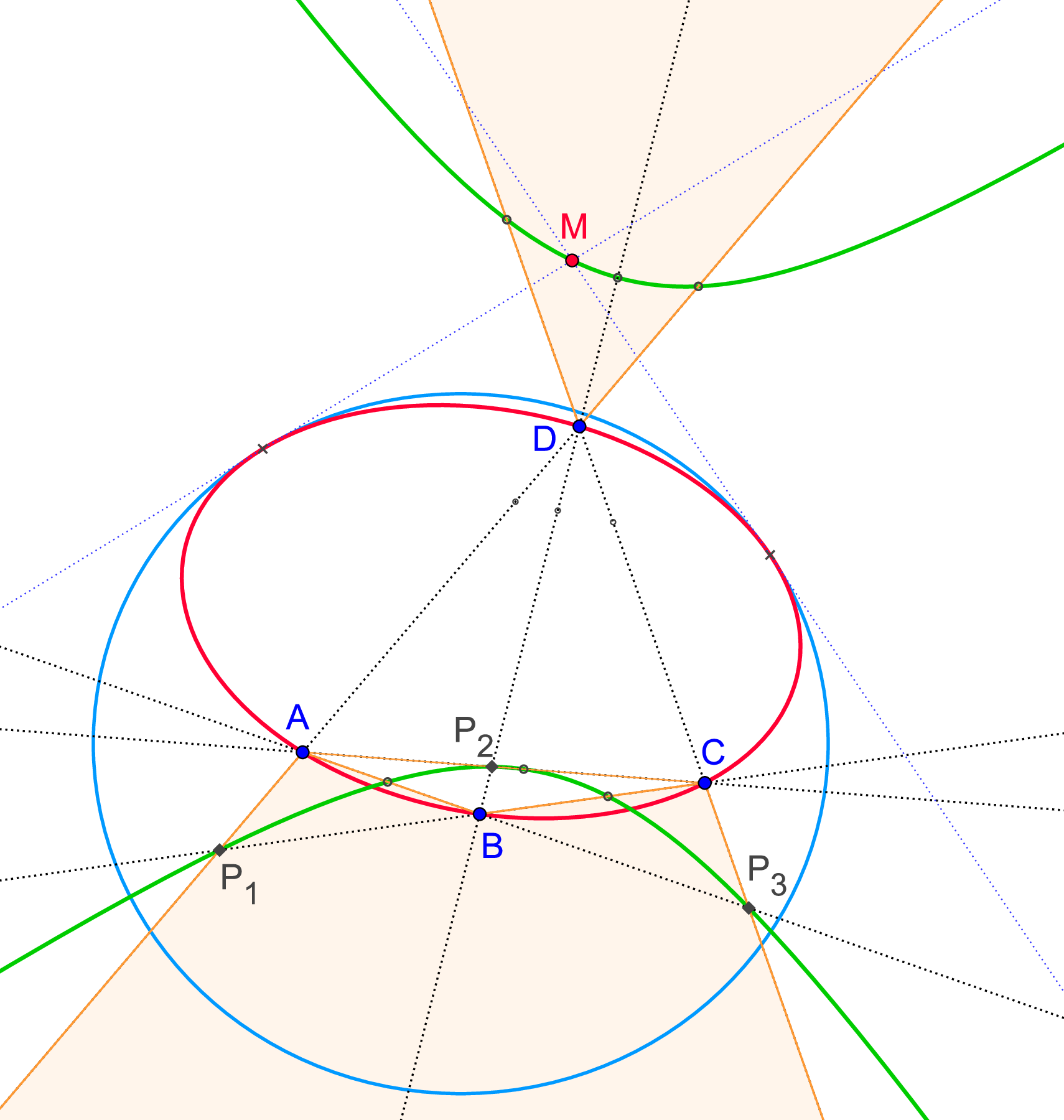}
\caption{A ten-point conic (green) in a hyperbolic plane. The blue circle is the absolute conic. The point M is the center of the red circle and a point on the conic. The segments connecting two vertices of the complete quadrangle do not have to lie inside the circle, nor does the circle center.
}
\includegraphics[height=7.6cm]{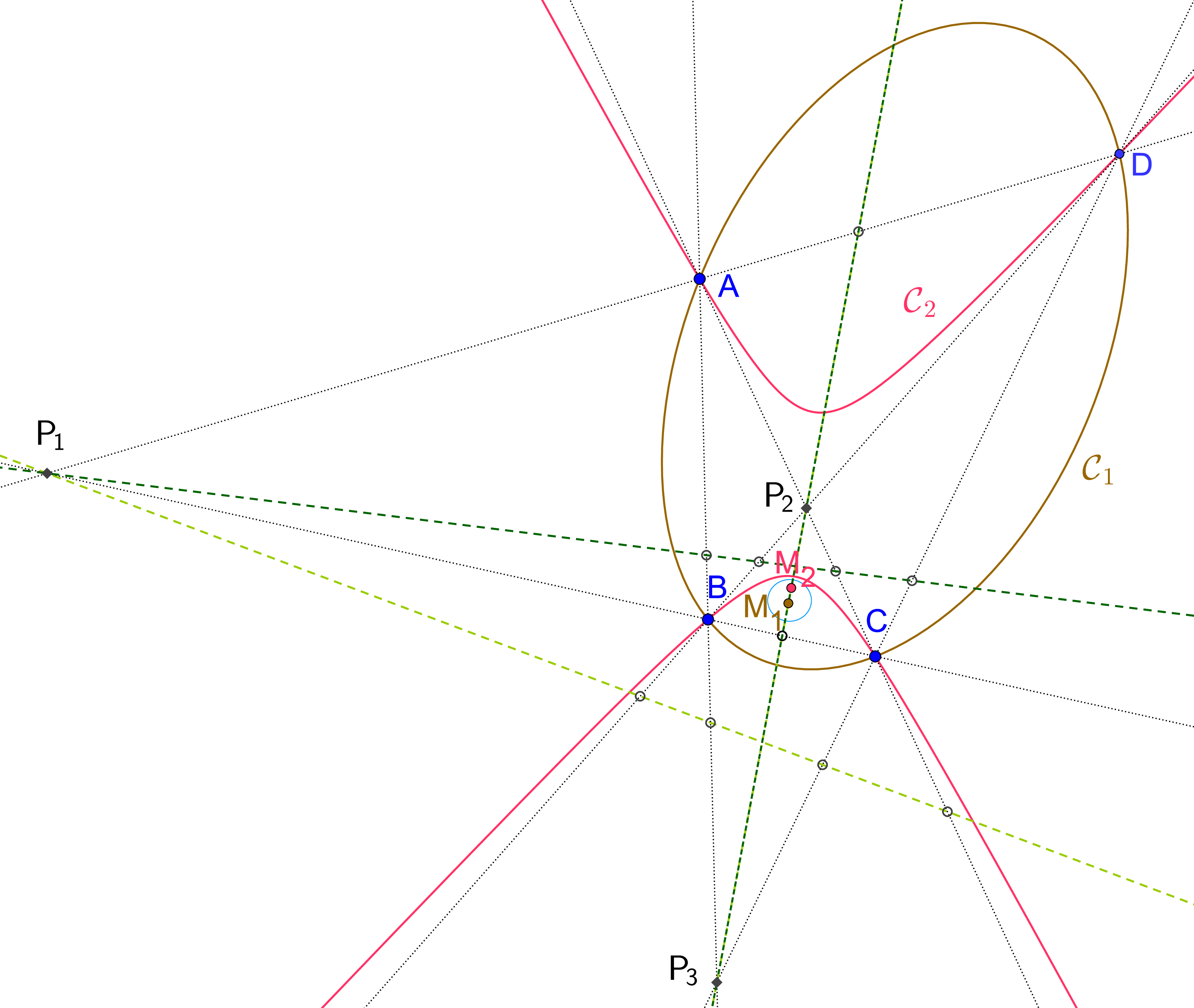}
\caption{Two ten-point conics (light green and dark green) of a tetragon $\{A,B,C,D\}$ in a hyperbolic plane.  Both conics are singular and share one line. This is the line through two distinct points which are the centers of two circles each of which passes through $\{A,B,C,D\}$. The light green conic "belongs" to the circle $\mathcal{C}_1$, the dark green conic to the circle $\mathcal{C}_2$. The small light blue circle is the absolute conic.}
\end{centering}
\vspace*{-1.5 mm}\end{figure}

What can we say about a nine-point conic in a regular Cayley-Klein plane? Let us find out.\vspace*{1 mm}\\
We call a complete graph (in the sence of graph-theory) with vertices in a projective plane a \textit{complete} \textit{tetragraph} if its vertices form a quadrangle and its edges are line-segments; and we call the edges of this graph the \textit{sides} of the complete tetragraph. \\
For a quadrangle $\mathbf{QA}=\{A,B,C,D\}$ in a projective plane there are 64 complete tetragraphs with these vertices $A,B,C,D$.\\
A \textit{six-point conic} of a quadrangle $\mathbf{QA}=\{A,B,C,D\}$ in an elliptic or in a hyperbolic plane is a conic which passes through all inner semi-midpoints of sides of a complete tetragraph of $\mathbf{QA}$. It is obvious that a quadrangle cannot have a six-point conic unless all its vertices are anisotropic. \\
A six-point conic of a complete quadrangle $\mathbf{QA}=\{A,B,C,D\}$ is called \textit{nine-point conic} of this complete quadrangle if it also goes through all three diagonal points of $\mathbf{QA}$.\\

In a regular Cayley-Klein plane, a six-point conic of a quadrangle $\mathbf{QA}=\{A,B,C,D\}$ need not be a nine-point conic. We give an example for the elliptic plane.\\
We may assume that the metric structure is given by the matrix {$\mathfrak{G} =\mathrm{diag}(1,1,1)$ with respect to the self-dual triangle $E_1E_2E_3$, $E_1{{\,{{\scriptstyle{:}}\!\!=}\,}}[1{:}0{:}0]$, $E_2{{\,{{\scriptstyle{:}}\!\!=}\,}}[0{:}1{:}0]$, $E_3{{\,{{\scriptstyle{:}}\!\!=}\,}}[0{:}0{:}1]$.\\ Put $A{{\,{{\scriptstyle{:}}\!\!=}\,}}[-3{\,:\,}0{\,:\,}4]$, $B{\;{\,{{\scriptstyle{:}}\!\!=}\,}\;}[0{\,:\,}3{\,:\,}4]$, $C{{\,{{\scriptstyle{:}}\!\!=}\,}}[75{\,:\,}{-24}{\,:\,}32]$ , $D{{\,{{\scriptstyle{:}}\!\!=}\,}}$  $[0{\,:\,}{-3}{\,:\,}4]$.\\ Then  $\mathfrak{M}= \left(\begin{array}{ccc} 
-43 & 16 & 6\\ 
16 & 75 & 6\\ 
6 & 6 &0
\end{array}\right)\;\;$}
is a (singular) matrix of a six-point conic of the quadrangle $\mathbf{QA}=\{A,B,C,D\}$, but not of a nine-point conic. This conic passes through the inner midpoints of $[A,B]_+, [A,C]_+, [A,D]_+, [B,C]_+, [B,D]_+, [C,D]_+ $, but not through any of the diagonal points $P_1=[{-50}{:}41{:}12], P_2=[0{:}{-2}{:}11], P_3=[{-6}{:}{-3}{:}4]$.\\

\noindent\textit{Theorem} 4. Let $\mathbf{QA}=\{A,B,C,D\}$ be a quadrangle in an elliptic or in a hyperbolic plane. We assume that $A,B,C,D$ are anisotropic points and $\{A\}{\,\cong\,}\{B\}{\,\cong\,}\{C\}{\,\cong\,}\{D\}$. 
Then there exists a nine-point conic for $\mathbf{QA}$ if and only if the points $A,B,C,D$ are concyclic.\\
When a circle passes through $A,B,C,D$, then the circle center is also a noteworthy point on the nine-point conic, so we may speak of a ten-point conic.\vspace*{0.7 mm}\\
\textit{Remark}. There can be two circles (not more) passing through these vertices. In this case, there are two complete tetragraphs belonging to $\mathbf{QA}$, each of them connected with a nine-point conic. There is a 1 to 1 relationship between the circles and the nine-point conics.  
\vspace*{1.0 mm}\\

\noindent \textit{Proof} of Theorem 4.
Let $A=[1{:}0{:}0], B=[0{:}1{:}0], C=[0{:}0{:}1]$ generate an elliptic or a hyperbolic plane with metric tensor $\mathfrak{G}$. 
Since $\{A\}{\,\cong}\,\{B\}{\,\cong\,}\{C\}{\,\cong\,}\{D\}$, we may assume that $\mathfrak{g}_{11}=\mathfrak{g}_{22}=\mathfrak{g}_{33}=1$, and for the point $D=[d_1{:}d_2{:}d_3]$ we have  $0<(d_1,d_2,d_3){\,\scriptscriptstyle{[\mathfrak{G}]}\,}(d_1,d_2,d_3){\,{=\!\!{\scriptstyle{:}}}\,}d^2$.\\
$A,B,C$ are  vertices of four triangles, 
$\Delta_0(A,B,C){\,{{\scriptstyle{:}}\!\!=}\,}\{[q_1{:}q_2{:}q_3]|q_1,q_2,q_3\geq 0\}$,\\
$\Delta_1(A,B,C)$ ${\,{{\scriptstyle{:}}\!\!=}\,}\{[q_1{:}q_2{:}q_3]$ $|q_1\geq 0,{q_2,q_3\leq 0}\}$, $\,\dots\,$.\\
The circumcircle of triangle $\Delta_0(A,B,C)$ is\vspace*{1 mm}\\
\centerline{$\mathcal{C}_0 =\{ [x_1{:}x_2{:}x_3]\;|\;(1-\mathfrak{g}_{12}) x_1 x_2 + (1-\mathfrak{g}_{13}) x_1 x_3 + (1-\mathfrak{g}_{23}) x_2 x_3 = 0\,\}$\;.}  \vspace*{-0.7 mm}\\ 
Thus,\; $\mathfrak{C}_0(A,B,C){\,{{\scriptstyle{:}}\!\!=}\,}\left(\begin{array}{ccc} 
0 & 1-\mathfrak{g}_{12} & 1-\mathfrak{g}_{13}\vspace*{0.5 mm} \\ 
1-\mathfrak{g}_{12} & 0 &1-\mathfrak{g}_{23}\vspace*{0.5 mm}\\ 
1-\mathfrak{g}_{13}&1-\mathfrak{g}_{23}&0
\end{array}\right)\;\;$ is a matrix of the circumcircle.\vspace*{0.9 mm}\\
Its center, the circumcenter of $\Delta_0$, is the point\vspace*{1.5 mm}\\
\noindent$O_0{\,{{\scriptstyle{:}}\!\!=}\,}[(1{-}\mathfrak{g}_{23})(1{+}\mathfrak{g}_{23}{-}\mathfrak{g}_{13}{-}\mathfrak{g}_{12}){:}
(1{-}\mathfrak{g}_{13})(1{+}\mathfrak{g}_{13}{-}\mathfrak{g}_{12}{-}\mathfrak{g}_{23}){:}(1{-}\mathfrak{g}_{12})(1{+}\mathfrak{g}_{12}{-}\mathfrak{g}_{23}{-}\mathfrak{g}_{13})]$.\vspace*{1.5 mm}\\
The circumcircle of $\Delta_1(A,B,C)$, $\mathcal{C}_1(A,B,C)$, has the center\vspace*{0.5 mm}\\
$O_1{\,{{\scriptstyle{:}}\!\!=}\,}[(1{-}\mathfrak{g}_{23})(1{+}\mathfrak{g}_{23}+\mathfrak{g}_{13}{+}\mathfrak{g}_{12}){:}(1{+}\mathfrak{g}_{13})(1{-}\mathfrak{g}_{13}{+}\mathfrak{g}_{12}{-}\mathfrak{g}_{23}){:}(1{+}\mathfrak{g}_{12})(1{-}\mathfrak{g}_{12}{+}\mathfrak{g}_{23}{-}\mathfrak{g}_{13})]$\\
\vspace*{0.7 mm}\\ and a matrix $\mathfrak{C}_1(A,B,C){\,{{\scriptstyle{:}}\!\!=}\,}\left(\begin{array}{ccc} 
0 & 1-\mathfrak{g}_{12} & 1-\mathfrak{g}_{13}\vspace*{0.5 mm} \\ 
1-\mathfrak{g}_{12} & 0 &1+\mathfrak{g}_{23}\vspace*{0.5 mm}\\ 
1-\mathfrak{g}_{13}&1+\mathfrak{g}_{23}&0
\end{array}\right),\;\;$ \vspace*{1.2 mm}\\
etc.\vspace*{0.8 mm}\\
The diagonal points of quadrangle $\{A,B,C,D\}$ are $P_1 = (A\vee D)\wedge(B\vee C) = [0{:}d_2{:}d_3]$, $P_2 = (A\vee C)\wedge(B\vee D) = [d_1{:}0{:}d_3]$, $P_3 = (C\vee D)\wedge(A\vee B) = [d_1{:}d_2{:}0]$.\vspace*{0 mm}\\

Assume now that
there is a conic passing through the inner midpoint $[1{:}1{:}0]$ of $[A,B]_+$, the inner midpoint $[0{:}1{:}1]$ of $[B,C]_+$ and through the points $P_1,P_2,P_3$. A matrix of this conic is \vspace*{1.9 mm}\\
\centerline{$ \mathfrak{M}_0 = \left(\begin{array}{ccc} 
2d_2 d_3 &-d_3(d_1+d_2)&-d_2(d_3+d_1)\vspace*{1 mm}\\ 
-d_3(d_1+d_2)&2d_3 d_1&-d_1(d_2+d_3)\vspace*{1 mm}\\ 
-d_2(d_3+d_1)&-d_1(d_2+d_3)&2d_1 d_2
\end{array}\right).$}\vspace*{1.9 mm}\\
This conic also passes through the inner midpoint $[1{:}0{:}1]$ of $[A,C]$, and it passes through
\centerline{- the inner midpoint $[d_1+d{:}d_2{:}d_3]$ of $[A,D]_+$ iff $(d_1+d_2+d_3-d)(d_3+d) = 0$,}
\centerline{- the inner midpoint $[d_1{:}d_2+d{:}d_3]$ of $[B,D]_+$ iff $(d_1+d_2+d_3-d)(d_2+d) = 0$,}
\centerline{- the inner midpoint $[d_1{:}d_2{:}d_3+d]$ of $[C,D]_+$ iff $(d_1+d_2+d_3-d)(d_1+d) = 0$,}
\centerline{- the outer midpoint $[d_1-d{:}d_2{:}d_3]$ of $[A,D]_+$ iff $(d_1+d_2+d_3+d)(d_3-d) = 0$,}
\centerline{- the outer midpoint $[d_1{:}d_2-d{:}d_3]$ of $[B,D]_+$ iff $(d_1+d_2+d_3+d)(d_2-d) = 0$,}
\centerline{- the outer midpoint $[d_1{:}d_2{:}d_3-d]$ of $[C,D]_+$ iff $(d_1+d_2+d_3+d)(d_1-d) = 0$.}\\
						
Thus, this conic is a nine-point conic of quadrangle $\{A,B,C,D\}$ precisely when $d^2\;=\,(d_1+d_2+d_3)^2.$  
We now show that this equation is a necessary and sufficient condition for the point $D$ to lie on the circumcircle of triangle $\Delta_0$:\\
{$ \hspace*{32.5mm} (d_1,d_2,d_3)\,{\scriptstyle{[\mathfrak{C}_0]}}(d_1,d_2,d_3) = 0$}\hspace*{14.9 mm}\\
\centerline{$\hspace*{-1.5 mm}\Leftrightarrow\;\;(d_1+d_2+d_3)^2 - (d_1,d_2,d_3){\scriptstyle{[\mathfrak{G}]}}(d_1,d_2,d_3) = 0$}\\
{$\hspace*{26.1mm}\Leftrightarrow\;\;(d_1+d_2+d_3)^2 - d^2 = 0$\vspace*{1.5 mm}}\\
\noindent\textit{Remark}. If $d=d_1+d_2+d_3$, $B$ is a point on the circumcircle of $\Delta_0(A,C,D)$. In the case of $d=\,-d_1-d_2-d_3$, $B$ lies on the circumcircle of $\Delta_3(A,C,D)$.\hspace*{0.5 mm}\\
 
It can be easily checked that if $d^2=(d_1+d_2+d_3)^2$, the circumcenter $O_0$  is also a point on the nine-point conic. \vspace*{1.0mm}\\

We now assume that $D$ is a point on the circumcircle of $\Delta_1(A,B,C)$. The matrix of the conic which passes through the outer midpoints of $[A,B]_+$ and $[A,C]_+$ and through the diagonal points $P_1,P_2,P_3$
is {$ \mathfrak{M}_1 = \left(\begin{array}{ccc} 
-2d_2 d_3 &d_3(d_1-d_2)&d_2(d_1-d_3)\vspace*{1 mm}\\ 
d_3(d_1-d_2)&2d_1 d_3&-d_1(d_2+d_3)\vspace*{1 mm}\\ 
d_2(d_3-d_1)&-d_1(d_2+d_3)&2d_1 d_2
\end{array}\right).$}\vspace*{0.9 mm}\\

This conic passes also through the inner midpoint of $[B,C]_+$. The outer midpoint of $[A,D]_+$ and the inner midpoints of $[B,D]_+$ and $[C,D]_+$ are also points on this conic iff $d+d_1-d_2-d_3=0$. The inner midpoint of $[A,D]_+$ and the outer midpoints of $[B,D]_+$ and $[C,D]_+$ lie on this conic iff $d-d_1+d_2+d_3=0$.
So this conic is a nine-point conic iff $(d+d_1-d_2-d_3)(d-d_1+d_2+d_3)=0$. But this equation is equivalent to the equation $(d_1,d_2,d_3){\scriptstyle{[\mathfrak{C}_1]}}(d_1,d_2,d_3)=0$, where $\mathfrak{C}_1$ is the matrix of the circumcircle of triangle $\Delta_1(A,B,C)$.\hspace*{1.5 mm}\\
If $d^2=(d_1-d_2-d_3)^2$, the circumcenter of $\Delta_1(A,B,C)$ is a tenth point on the conic.\\

%We will spare ourselves the discussion of the remaining two cases ($D$ on the circumcircle of $\Delta_i(A,B,C), i=2,3.$  $\;\Box$

In the proof so far, two cases have been examined. Any two points of the quadrangle $\{A,B,C,D \}$ were connected by a segment to get a complete tetragraph  ; in the first case we chose the line segments $[A,B]_+ , [A,C]_+ , [A,D]_+ , [B,C]_+ , [B,D]_+ , [C,D]_+$, in the second the line segments 
$[A,B]_+ , [A,C]_+ , [A,D]_- , [B,C]_+ , [B,D]_- , [C,D]_-$. We do not give an explicit proof for the other two cases, since the proof is done quite analogously.$\;\;\;\Box$\\

\noindent \textit{Addendum} to Theorem 4. If the vertices $A, B, C, D$  of a quadrangle $\mathbf{QA}$ lie on a circle with center $O$, then there exist a nine-point conic and a complete tetragraph of $\mathbf{QA}$ such that the  conic runs through the inner midpoints of the sides of the tetragraph. It can be easily checked that all the outer midpoints lie collinearly on a line $\mathcal{L}=\mathcal{L}(\mathbf{QA},\mathcal{C})$, which is the dual of $O$. 
Animate a point $P$ on this line; the trace of the intersection  $P'$ of the dual of $P$ with polar of the circle is a conic, cf. \cite[p. 342]{Pe}. This conic must be the nine-point conic, because if $P$ is an outer midpoint of a side of the complete tetragraph, both the dual line of $P$ and the polar of $P$ with respect to the circle pass through the inner midpoint.

The nine-point conic of $(\mathbf{QA},\mathcal{C})$ intersects the line $\mathcal{L}$ at two antipodal points $Q_1$ and $Q_2$.\\
If this nine-point conic is nonsingular, any conic passing through the points $A, B, C, D$ touches $\mathcal{L}$ if and only if it passes through $Q_1$ or $Q_2$, see \cite{Al}.\\
If the nine-point conic is singular, then there are two distinct circles $\mathcal{C}_1$ and $\mathcal{C}_2$ through $A, B, C, D$. The two lines $\mathcal{L}_i = \mathcal{L}_i(\mathbf{QA},\mathcal{C}_i), i=1,2 ,$ meet at a point $Q_1$ which is a center of a singular conic through $A, B, C, D$; and there are two nonsingular conics passing through $A, B, C, D$, one touching $\mathcal{L}_1$ at the antipodal point of $Q_1$ on $\mathcal{L}_1$ and the other touching $\mathcal{L}_2$ at the antipodal point of $Q_1$ on $\mathcal{L}_2$.\\ 
In any case, these two points $Q_1$, $Q_2$ are two more noteworthy points on the nine-point conic.

If $A,B,C,D$ are anisotropic points in a hyperbolic plane and $\{A\}{\,\cong\,}\{B\}{\,\cong\,}\{C\}{\,\ncong\,}\{D\}$, then there cannot be any circle passing through all four points. But the quadrangle $\mathbf{QA}=\{A,B,C,D\}$ can have a nine-point conic. Here is an example. \\
We may assume that the metric structure in a hyperbolic plane is given by the matrix {$\mathfrak{G} =\mathrm{diag}(1,1,-1)$ with respect to the triple $(E_1,E_2,E_3)$ of points $E_1{{{{\scriptstyle{:}}\!\!=}\,}}[1{:}0{:}0]$, $E_2{{{{\scriptstyle{:}}\!\!=}\,}}[0{:}1{:}0]$, $E_3{{{{\scriptstyle{:}}\!\!=}\,}}[0{:}0{:}1]$. Put $A:=[{-1}{\;:\,}{-1}{\;:\,}2], B:=[0{\;:\,}1{\;:\,}4],$ $ C:=[1{\;:\,}{-1}{\;:\,}2], D:=[6 \sqrt{37607 - 6866 \sqrt{30}} - 150 \sqrt{30} + 820{\,:\,} - 3 \sqrt{37607 - 6866 \sqrt{30}} + 75 \sqrt{30} - 413{\,:\,} 4]$.\\The outer semi-midpoints of the segments $[A,B]_+$ , $[A,C]_+$ , $[A,D]_+$ , $[B,C]_+$ , $[B,D]_+ ,$ $[C,D]_+$ lie collinearly on the line $[0{\,:\;}{-19}-3\sqrt{30}{\,:\;}14] \vee [1{\,:\;}0{\,:\;}0]$, so there is a conic passing through the inner semi-midpoints of these segments, see  Fig. 10. 
\begin{figure}[!htbp]
\begin{centering}
\includegraphics[height=5.7cm]{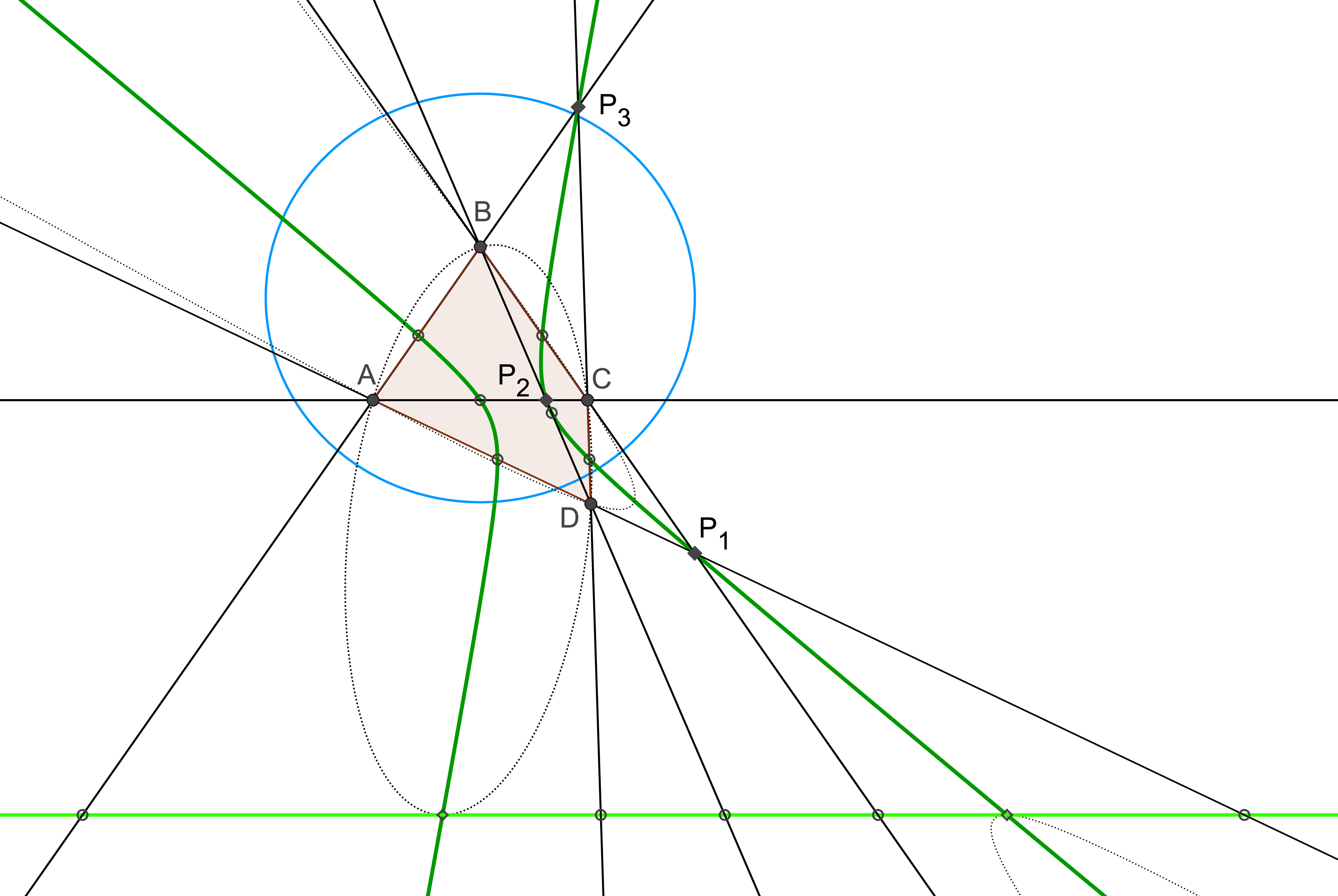}
\caption{}
\end{centering}
\vspace*{1 mm}\end{figure}
%\noindent\textit{Definition}. Symmetry points and symmetry axes of real conics.\\
%An anisotropic point $P = [p_1{:}p_2{:}p_3]$ is a symmetry point of a nonempty real conic $\mathcal{C}$ if for every point
%$Q = [q_1{:}q_2{:}q_3]$ on this conic the mirror image of $Q$ in $P$ is also a point on this
%conic. A line $\mathcal{L}$ is a symmetry axis of $\mathcal{C}$ if its dual point is a symmetry point. 
%Here are four examples.
%\hspace*{5mm}If the conic is the union of two different lines which meet at an anisotropic point $P$, then the point $P$ and the duals of the two angle bisectors are the symmetry points. The point $P$ is regarded as the center of this conic. \\
%\hspace*{5mm}If the conic is a nonisotropic line $\mathcal{L}$ with multplicity 2, the dual of $\mathcal{L}$ and each point $P$ on $\mathcal{L}$ are symmetry points.\\
%\hspace*{5mm}The symmetry points of a proper circle with center $O$ are the center $O$ and all points on the dual line $P^\delta$.\\ 
%\hspace*{5mm}A nonsingular conic which does not touch the absolute conic has three symmetry points, one lies inside the conic and is regarded as its center, the other two lie outside on the dual line of the center.\\
%Remark. A nonsingular conic may have up to four isotropic points.\\

\noindent \textit{Theorem} 5.  We assume that $\mathbf{QA}$ is a quadrangle with anisotropic vertices $A, B, C, D$ on a circle $\mathcal{C}$. If $\mathcal{K}$ is any conic passing through $A, B, C, D$, then either $\mathcal{K}=\mathcal{C}$ or all its centers are points on the nine-point conic of $\mathbf{QA}$.\\
\textit{Proof}. We only provide a proof for the case that $D=[d_1{:}d_2{:}d_3]$ is a point on $\mathcal{C}=\mathcal{C}_0(A,B,C)$. (For the remaining three cases the proof is quite analogous.)\\
If $D\in \mathcal{C}_0(A,B,C)$, we can find a real number $t$ such that $d_1=t\, (t\,(\mathfrak{g}_{12}-1)+\mathfrak{g}_{23}-1), d_2=t\,(1-\mathfrak{g}_{13}), d_3=t\,(\mathfrak{g}_{12}-1)+\mathfrak{g}_{23}-1$.
Let $M=[m_1{:}m_2{:}m_3]$ be an anisotropic point and let $A',B',C',D'$ be the mirror images of $A, B, C, D$ under a reflection in $M$.\\
The conic through $A, B, C, D, A'$ is $\mathcal{F}=\{ [x_1{:}x_2{:}x_3]\,|\,f_{\!12}x_1 x_2+ f_{\!13}x_1 x_3+f_{\!23}x_2 x_3\,\}$ with\\ 
$f_{\!12}{\,{{\scriptstyle{:}}\!\!=}\,}m_3 (2 t m_2 ((\mathfrak{g}_{12} - 1) (m_1 + \mathfrak{g}_{12} m_2 + \mathfrak{g}_{13} m_3)) + 2 \mathfrak{g}_{12} m_2^2 (\mathfrak{g}_{23} - 1) +\\
\hspace*{10mm}+ \mathfrak{g}_{13} (m_1^2 - m_2^2 - 2 m_2 m_3 - m_3^2) + 2 \mathfrak{g}_{23} m_2 (m_1 + m_3) - m_1^2 - 2 m_1 m_2 + m_2^2 + m_3^2),$\\
$f_{\!13}{\,{{\scriptstyle{:}}\!\!=}\,}m_2 (\mathfrak{g}_{13} - 1) (2 t m_3 (\mathfrak{g}_{13} m_2 + \mathfrak{g}_{13} m_3 + m_1) + 2 \mathfrak{g}_{23} m_2 m_3 - m_1^2 + m_2^2 + m_3^2)$,\\
$f_{\!23}{\,{{\scriptstyle{:}}\!\!=}\,}(2 \mathfrak{g}_{23} m_2 m_3 - m_1^2 + m_2^2 + m_3^2) (t\,(\mathfrak{g}_{12} m_2 + \mathfrak{g}_{13} m_3 - m_2 - m_3) + m_2 (\mathfrak{g}_{23} - 1))$.\vspace*{1.5mm}\\
Using a computer algebra system (CAS) we find a condition for $B',C',D'$ to also lie on $\mathcal{F}$: 
\centerline{$(m_1+m_2+m_3)\,f\!(m_1,m_2,m_3)= 0$, with}
$f\!(m_1,m_2,m_3){\,{{\scriptstyle{:}}\!\!=}\,}t\,\big(\mathfrak{g}_{12} m_2 (m_1{-}m_2{+}m_3)-\mathfrak{g}_{13} m_3 (m_1{+}m_2{-}m_3) + (m_1{-}m_2{-}m_3) (m_3{-}m_2)\big)+\\
\hspace*{5mm}+\mathfrak{g}_{13} m_1 (m_1{-}m_2-m_3)+\mathfrak{g}_{23} m_2 (m_1-m_2{+}m_3)-m_1(m_1{-}m_3)+m_2(m_2{-}m_3)$.\\
If $m_1+m_2+m_3=0$, $M$ is a point on the dual line of the center of $\mathcal{C}$ and therefore $\mathcal{F}=\mathcal{C}$.\\
If $f\!(m_1,m_2,m_3)= 0$, $M$ is a point on the nine-point conic. $\;\;\;\Box$\\

\begin{figure}[!tbph]
\begin{centering}
\includegraphics[height=8.6cm]{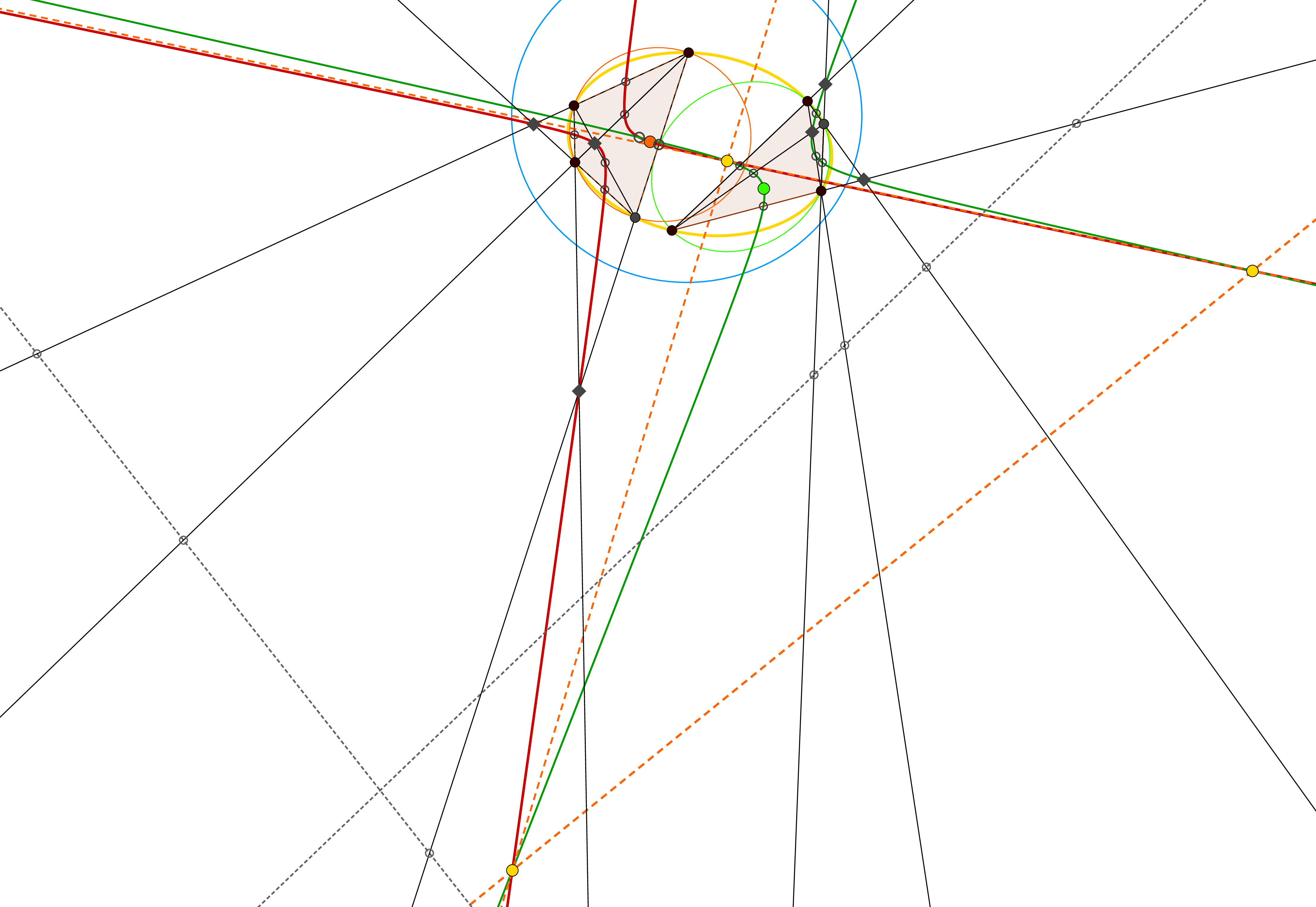}
\caption{This figure shows how the symmetry points of an ellipse in the hyperbolic plane can be constructed by the intersection of two nine-point conics. The yellow points are the symmetry points of the yellow ellipse. It is possible to create a tool in GeoGebra which, after entering a conic section, outputs its symmetry points. 
}
\end{centering}
\vspace*{1 mm}\end{figure}

\noindent \textit{Theorem} 6. The center of a circle $\mathcal{C}$ touching the sidelines of a tetragon  $\mathbf{T} = ABCD$ lies on a conic which passes through the  inner and outer midpoints of the segments $[A{\,,\,}C]_+ , [B{\,,\,}D]_+$, $[P_1{\,,\,}P_3]_+$, where $P_1, P_2, P_3$ are the diagonal points and $A{\vee}C$, $B{\vee}D$, $P_1{\vee}P_3$ the diagonals of $\mathbf{T}$. See Figure 12.\vspace*{0.0 mm}\\

\begin{figure}[!bhpt]
\begin{centering}
\includegraphics[height=10.2cm]{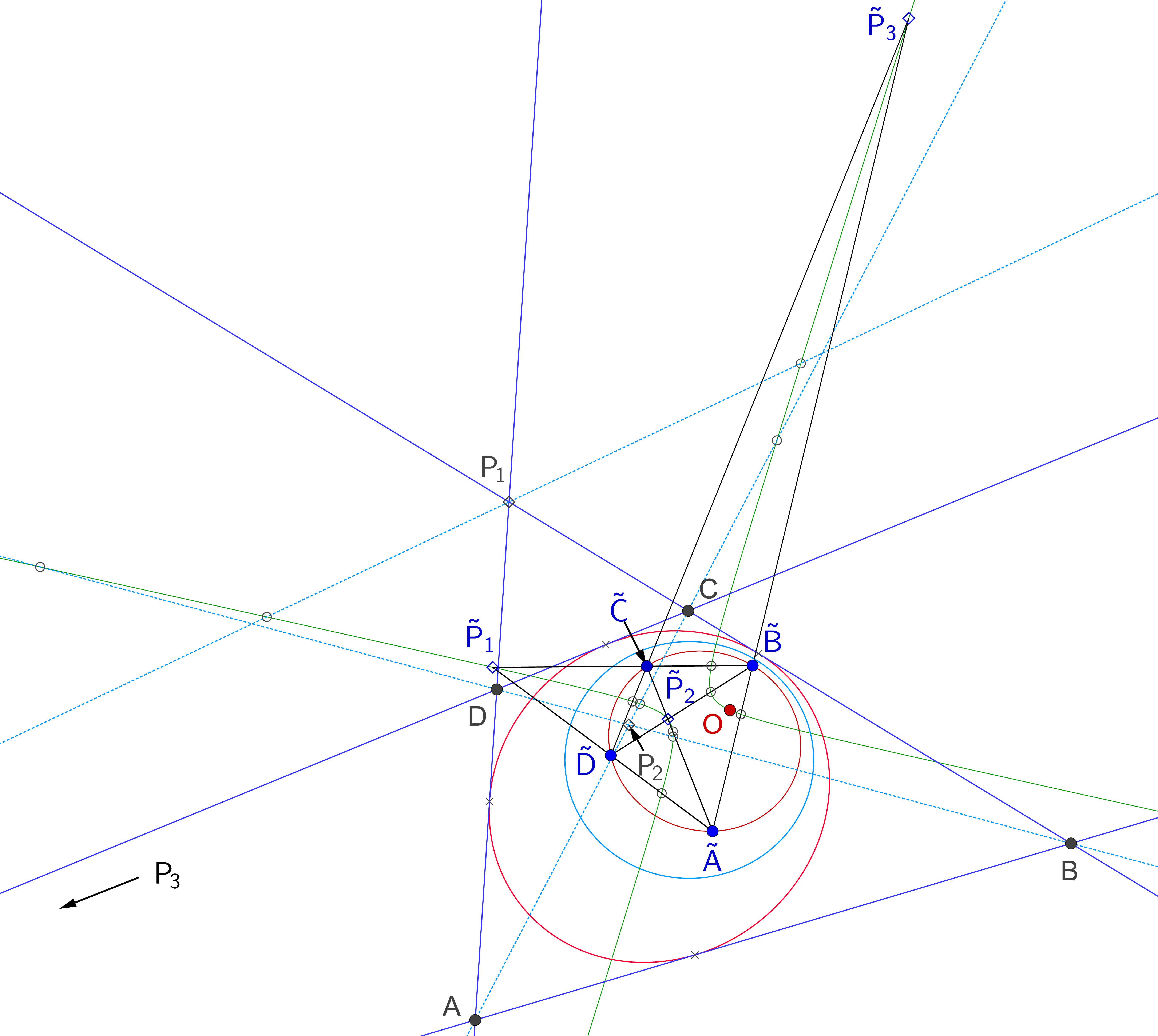}
\caption{}
\end{centering}
\vspace*{1 mm}\end{figure}

\noindent\textit{Proof}:
If $\mathbf{T} = (A,B,C,D)$ is a tetragon and $\mathcal{C}$ a circle tangent to its sidelines, then the points $\tilde{A}{\,{{\scriptstyle{:}}\!\!=}\,}(A\vee B)^{\mathfrak{G}}$, $\tilde{B}{\,{{\scriptstyle{:}}\!\!=}\,}(B\vee C)^{\mathfrak{G}}$, $\tilde{C}{\,{{\scriptstyle{:}}\!\!=}\,}(C\vee D)^{\mathfrak{G}}$ and $\tilde{D}{\,{{\scriptstyle{:}}\!\!=}\,}(D\vee A)^{\mathfrak{G}}$ are points on a circle $\tilde{\mathcal{C}}$ that is concentric to $\mathcal{C}$, see Figure 11.\vspace*{0.5 mm}

We take $\tilde{A}, \tilde{B}, \tilde{C}$ as a generating system for the plane. It is not difficult to determine a number $i \in \{0,1,2,3\}$ such that $\tilde{{\mathcal{C}}}_i$ is a circumcircle of a  triangle $\Delta_i(\tilde{A}, \tilde{B}, \tilde{C})$ which passes through $\tilde{D}$. A matrix $\tilde{\mathfrak{K}}_i$ can be set up for the corresponding nine-point conic $\tilde{\mathcal{K}}_i$ of the quadrangle $\mathbf{QA}=\{\tilde{A}, \tilde{B}, \tilde{C},\tilde{D}\}$}.\vspace*{0.5 mm}

To prove the theorem, it is sufficient to show that $\tilde{\mathcal{K}}_i$ also passes through the inner and outer midpoints of the segments $[A,C]_+, [B,D]_+, [P_1,P_3]_+$. We will carry out this proof for the case that $i=0$; the proof for the other three cases is analogous.\vspace*{1.3 mm}\\
\noindent Since $\tilde{D}=[\tilde{d}_1{:}\tilde{d}_2{:}\tilde{d}_3]$ is a point on $\tilde{\mathcal{C}}_0$,\vspace*{1.3 mm}\\ 
\centerline{$\tilde{D} = [\tilde{d}_1(\mathfrak{g}_{13} \tilde{d}_1+\mathfrak{g}_{23}-\tilde{d}_1-1){\,:\,}\mathfrak{g}_{13} \tilde{d}_1+\mathfrak{g}_{23}-\tilde{d}_1-1{\,:\,}\tilde{d}_1 (1-\mathfrak{g}_{12})]$.}
\vspace*{-0.6 mm}\\
A suitable matrix for $\tilde{\mathcal{K}}_0$ is\vspace*{1.3 mm}\\ $\tilde{\mathfrak{K}}_0 =
{ \left(\begin{array}{ccc} 
2 (1{-\,}\mathfrak{g}_{12})&(\mathfrak{g}_{12}{-1})(d_1{+}1)&\mathfrak{g}_{12}{-\,}\mathfrak{g}_{13} d_1{-\,}\mathfrak{g}_{23}{+\,}d_1\vspace*{1 mm}\\ 
(\mathfrak{g}_{12}{-}1)(d_1{+}1)&2d_1(1{-\,}\mathfrak{g}_{12})&\mathfrak{g}_{12}d_1{-\,}\mathfrak{g}_{13} d_1{-\,}\mathfrak{g}_{23}{+}1\vspace*{1 mm}\\ 
\mathfrak{g}_{12}{-\,}\mathfrak{g}_{13} d_1{-\,}\mathfrak{g}_{23}{+\,}d_1&\mathfrak{g}_{12} d_1{-\,}\mathfrak{g}_{13} d_1{-\,}\mathfrak{g}_{23}{+}1&2(\mathfrak{g}_{13} d_1{+\,}\mathfrak{g}_{23}{-\,}d_1{-}1)
\end{array}\right).}$\vspace*{1.9 mm}\\

We calculate the barycentric coordinates of $A, B, C, D, P_1, P_3$ with respect to $(\tilde{A},\tilde{B},\tilde{C})$:
Put\\
${\mathbf{a}{\,{{\scriptstyle{:}}\!\!=}\,}\big((\mathfrak{g}_{23} - 1) (\mathfrak{g}_{12} \mathfrak{g}_{23} - \mathfrak{g}_{13}) {\;-\;} d_1 (\mathfrak{g}_{12}^2{\;-\;} \mathfrak{g}_{12} (2 \mathfrak{g}_{13} \mathfrak{g}_{23} - \mathfrak{g}_{23} + 1) {\;+\;} \mathfrak{g}_{13} (\mathfrak{g}_{13} + \mathfrak{g}_{23} - 1)),\;}\;\\d_1 (1{-}\mathfrak{g}_{13}) (\mathfrak{g}_{12} - \mathfrak{g}_{13} + \mathfrak{g}_{23} - 1) + (\mathfrak{g}_{23}{-}1) (\mathfrak{g}_{12} \mathfrak{g}_{13} - \mathfrak{g}_{23}), (\mathfrak{g}_{12}{-}1) (d_1 (\mathfrak{g}_{12} - \mathfrak{g}_{13} - \mathfrak{g}_{23} + 1) + (1 - \mathfrak{g}_{23}) (\mathfrak{g}_{12} + 1))\big)\,,$\\
$\mathbf{b}{\,{{\scriptstyle{:}}\!\!=}\,}\big(\mathfrak{g}_{13}-\mathfrak{g}_{12} \mathfrak{g}_{23}, \mathfrak{g}_{23}-\mathfrak{g}_{12} \mathfrak{g}_{13}, \mathfrak{g}_{12}^2-1\big)\,,$\\
$\mathbf{c}{\,{{\scriptstyle{:}}\!\!=}\,}\big(\mathfrak{g}_{23}^2-1, \mathfrak{g}_{12}-\mathfrak{g}_{13} \mathfrak{g}_{23}, \mathfrak{g}_{13}-\mathfrak{g}_{12} \mathfrak{g}_{23}\big)\,,$\\
$\mathbf{d}{\,{{\scriptstyle{:}}\!\!=}\,}\big(\mathfrak{g}_{12} d_1 - \mathfrak{g}_{13} \mathfrak{g}_{23} d_1 - \mathfrak{g}_{23}^2 + 1, -\mathfrak{g}_{12} + \mathfrak{g}_{13}^2 d_1 + \mathfrak{g}_{13} \mathfrak{g}_{23} - d_1, \mathfrak{g}_{12} (\mathfrak{g}_{23} - \mathfrak{g}_{13} d_1) - \mathfrak{g}_{13} + \mathfrak{g}_{23} d_1\big)\,,$\\
$\mathbf{p}_1{\,{{\scriptstyle{:}}\!\!=}\,}\big(d_1(\mathfrak{g}_{13}-1) (\mathfrak{g}_{12} \mathfrak{g}_{23}-\mathfrak{g}_{13})+(1-\mathfrak{g}_{23})(\mathfrak{g}_{12}+\mathfrak{g}_{13}-\mathfrak{g}_{23}-1), d_1(\mathfrak{g}_{13}-1) (\mathfrak{g}_{12} \mathfrak{g}_{13}-\mathfrak{g}_{23})-\mathfrak{g}_{12}^2+\mathfrak{g}_{12}(\mathfrak{g}_{13}(2 \mathfrak{g}_{23}-1)+1)-\mathfrak{g}_{23}(\mathfrak{g}_{13}+\mathfrak{g}_{23}-1), (1-\mathfrak{g}_{12})(d_1(\mathfrak{g}_{12}+1)(\mathfrak{g}_{13}-1)-\mathfrak{g}_{12}+\mathfrak{g}_{13}+\mathfrak{g}_{23}-1)\big)\,,$\\
$\mathbf{p}_3{\,{{\scriptstyle{:}}\!\!=}\,}\big(\mathfrak{g}_{12} - \mathfrak{g}_{13} \mathfrak{g}_{23},\mathfrak{g}_{13}^2-1, \mathfrak{g}_{23} - \mathfrak{g}_{12} \mathfrak{g}_{13}\big)\,.$\vspace*{1.5 mm}\\
There exist nonzero real numbers $a,b,c,d,p_1,p_3$ such that\vspace*{0.7 mm}\\
\centerline {$ A^\circ = \mathbf{a}/a{\,,\,}B^\circ = \mathbf{b}/{b}{\,,\,}C^\circ = {\mathbf{c}}/{c}{\,,\,}D^\circ = {\mathbf{d}}/{d}{\,,\,}P^\circ_1 = {\mathbf{p}_1}/{p_1}{\,,\,}P^\circ_3 = {\mathbf{p}_3}/{p_3}.\;\;\;\;$}\vspace*{1 mm}\\
We get (with the help of a computer algebra system) \vspace*{1.5 mm}\\
{\noindent \hspace*{5 mm}{$(a\,c)^2\,|(A^\circ \pm C^\circ){\,\scriptstyle{[\tilde{\mathfrak{K}}_0]}\,}(A^\circ \pm C^\circ)|$\vspace*{0.6 mm}\vspace*{0.5 mm}\\$=\; |(c\,  \mathbf{a} \pm  a\, \mathbf{c}){\,\scriptstyle{[\tilde{\mathfrak{K}}_0]}\,}(c\, \mathbf{a} \pm  a\, \mathbf{c})|$}\vspace*{1.5 mm}\\
\noindent \hspace*{0.2 mm}= $\,\big|\,2 d_1(\mathfrak{g}_{12} d_1-\mathfrak{g}_{13} d_1-\mathfrak{g}_{23}+1)(2 \mathfrak{g}_{12} \mathfrak{g}_{13} \mathfrak{g}_{23}-\mathfrak{g}_{12}^2-\mathfrak{g}_{13}^2-\mathfrak{g}_{23}^2+1)\,\cdot$ \\
\noindent \hspace*{6 mm}$\;\cdot\,\big(a^2\,(\mathfrak{g}_{23}+1)$\\
	\noindent \hspace*{10 mm}${+}\;
	c^2\;(1-\mathfrak{g}_{12})(\mathfrak{g}_{12}(\mathfrak{g}_{23}{-}2d_1{-}1)+2\mathfrak{g}_{13} d_1(d_1{+}1)+\mathfrak{g}_{23}(2d_1{+}1)-2d_1^2-2d_1-1)\big)  \big|$\\
\noindent \hspace*{0 mm}= $\displaystyle\,\big|\frac{2 d_1(\mathfrak{g}_{12} d_1-\mathfrak{g}_{13} d_1-\mathfrak{g}_{23}+1)(a^2 c^2 - c^2 a^2)}{(\mathfrak{g}_{23}-1)}\big|\;=\;0\,,$\vspace*{2.5 mm}\\
\noindent \hspace*{5 mm}{$(b\,d)^2\,|(B^\circ \pm D^\circ){\,\scriptstyle{[\tilde{\mathfrak{K}}_0]}\,}(B^\circ \pm D^\circ)|\vspace*{1 mm}\\ 
=\; |(d\,  \mathbf{b} \pm  b\, \mathbf{d}){\,\scriptstyle{[\tilde{\mathfrak{K}}_0]}\,}(d\,  \mathbf{b} \pm  b\, \mathbf{d})|$}\vspace*{0.6 mm}\\
\noindent \hspace*{0.2 mm}= $\,|\,2 d_1(d_1{+}1)(2 \mathfrak{g}_{12} \mathfrak{g}_{13} \mathfrak{g}_{23}-\mathfrak{g}_{12}^2{-}\,\mathfrak{g}_{13}^2{-\,}\mathfrak{g}_{23}^2{+}1)\,\cdot \vspace*{0.3 mm}$\\
\noindent \hspace*{6 mm}$\;\cdot\,\big(b^2\,(2 \mathfrak{g}_{12} d_1 - \mathfrak{g}_{13}^2 d_1^2-2 \mathfrak{g}_{13} \mathfrak{g}_{23} d_1{-\,}\mathfrak{g}_{23}^2{+\,}d_1^2{+\,}1) + d^2\,(\mathfrak{g}_{12}^2{-\,}1)  \big)|$\, \\
\noindent $\hspace*{0.2 mm}=\; |\,2 d_1(d_1{+}1)(b^2 d^2 - d^2 b^2)\,|\;=\;0\,,$\vspace*{2.5 mm}\\
\noindent \hspace*{5 mm}{$(p_1\,p_3)^2\,|(P_1^\circ \pm P_3^\circ){\,\scriptstyle{[\tilde{\mathfrak{K}}_0]}\,}(P_1^\circ \pm P_3^\circ)|\vspace*{1 mm}\\ 
=\; |\,(p_3\,  \mathbf{p}_1 \pm  p_1\, \mathbf{p}_3){\,\scriptstyle{[\tilde{\mathfrak{K}}_0]}\,}(p_3\,  \mathbf{p}_1 \pm  p_1\, \mathbf{p}_3)\,|$}\vspace*{0.6 mm}\\
\noindent \hspace*{0.2 mm}= $|\,2 d_1(\mathfrak{g}_{12} d_1-\mathfrak{g}_{13} d_1-\mathfrak{g}_{23}+1)(2 \mathfrak{g}_{12} \mathfrak{g}_{13} \mathfrak{g}_{23}-\mathfrak{g}_{12}^2-\mathfrak{g}_{13}^2-\mathfrak{g}_{23}^2+1)\,\cdot$ \\
\noindent \hspace*{6 mm}$\;\cdot\, (p_1^2\;((\mathfrak{g}_{12} - 1) (d_1^2 (\mathfrak{g}_{12}+1)(\mathfrak{g}_{13}-1)-2d_1(\mathfrak{g}_{12}-\mathfrak{g}_{13}-\mathfrak{g}_{23}+1)+2(\mathfrak{g}_{23} - 1)))\;-\;p_2^2\;(\mathfrak{g}_{13}+1)) |$\vspace*{1 mm}\\
\noindent \hspace*{0.2 mm}= $|\,2 d_1(\mathfrak{g}_{12} d_1-\mathfrak{g}_{13} d_1-\mathfrak{g}_{23}+1)(2 \mathfrak{g}_{12} \mathfrak{g}_{13} \mathfrak{g}_{23}-\mathfrak{g}_{12}^2-\mathfrak{g}_{13}^2-\mathfrak{g}_{23}^2+1)(p_1^2 p_3^2 - p_3^2 p_1^2)\,|\;=\;0\,.\;\;\Box$\vspace*{2 mm}\\
  \noindent\textit{Corollary}: In an elliptic and in a hyperbolic plane, each cyclic quadrilateral can be assigned a conic that passes through sixteen prominent points.\\

\noindent \textit{Theorem} 7.	Let $\mathcal{K}$ be a conic which touches the sides of a tangent quadrilateral $\mathbf{T} = ABCD$. If $\mathcal{K}$ is not a circle, then all symmetry points of $\mathcal{K}$ lie on a conic which passes through the inner and outer midpoints of the segments $[A{\,,\,}C]_{+} , [B{\,,\,}D]_{+} , [P_1{\,,\,}P_3]_{+}  $, where $P_1, P_2, P_3$ are the diagonal points and $A{\vee}C$, $B{\vee}D$, $P_1{\vee}P_3$ the diagonals of $\mathbf{T}$.\\
\textit{Proof}: The symmetry points of a conic agree with the symmetry points of its dual. So this theorem is a consequence of Theorem 5, Theorem 6 and the proof of Theorem 6. $\;\;\;\;\;\;\;\;\Box$

\section*{The Shatunov-Togarev line of a tetragon}

\begin{figure}[!htbp]
\begin{centering}
\includegraphics[height=5.8cm]{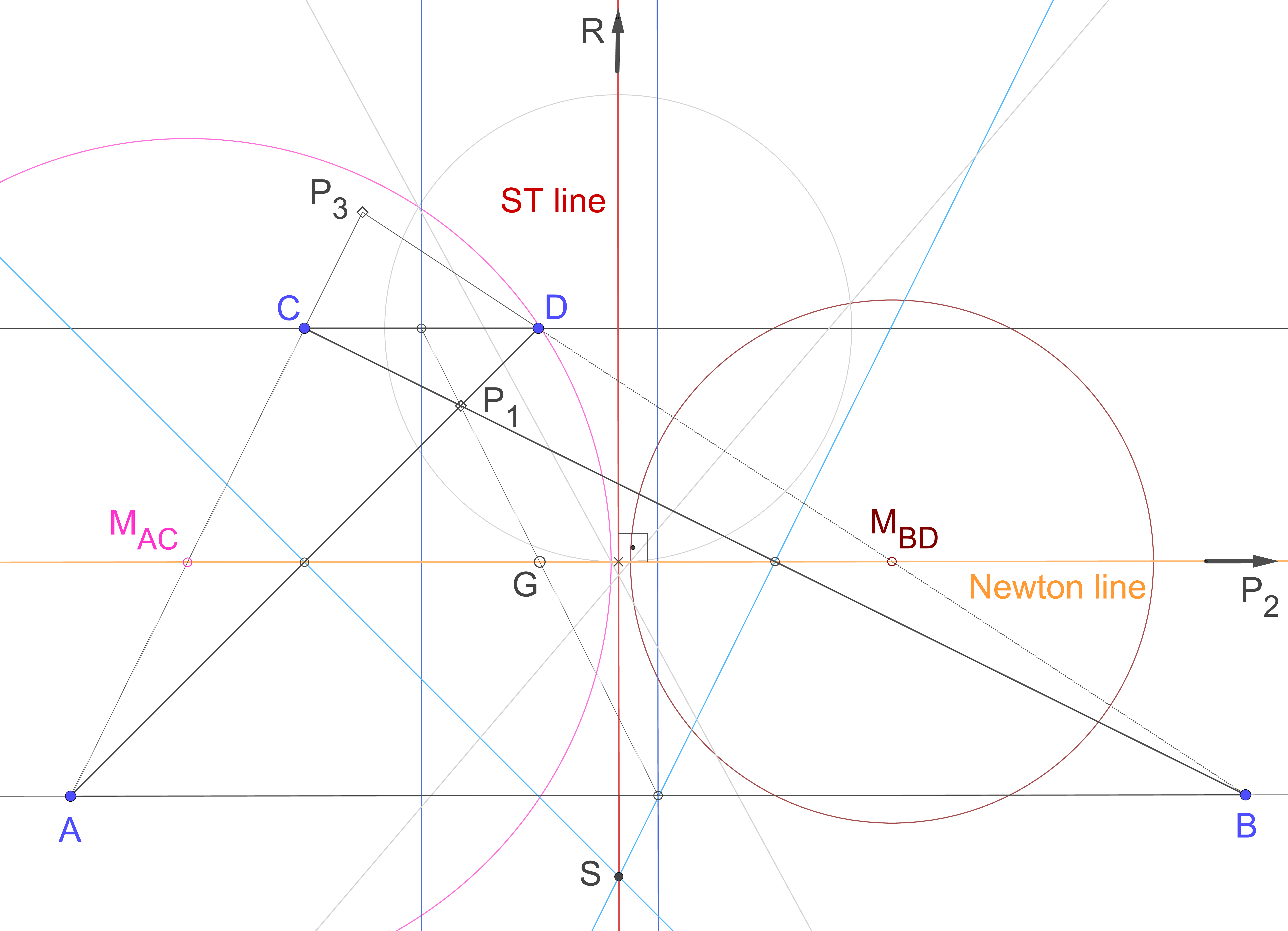}
\caption{Shatunov-Tokarev line and Newton line of a self-inter\-secting trapezoid in a Euclidean plane.}
\end{centering}
\includegraphics[height=5.8cm]{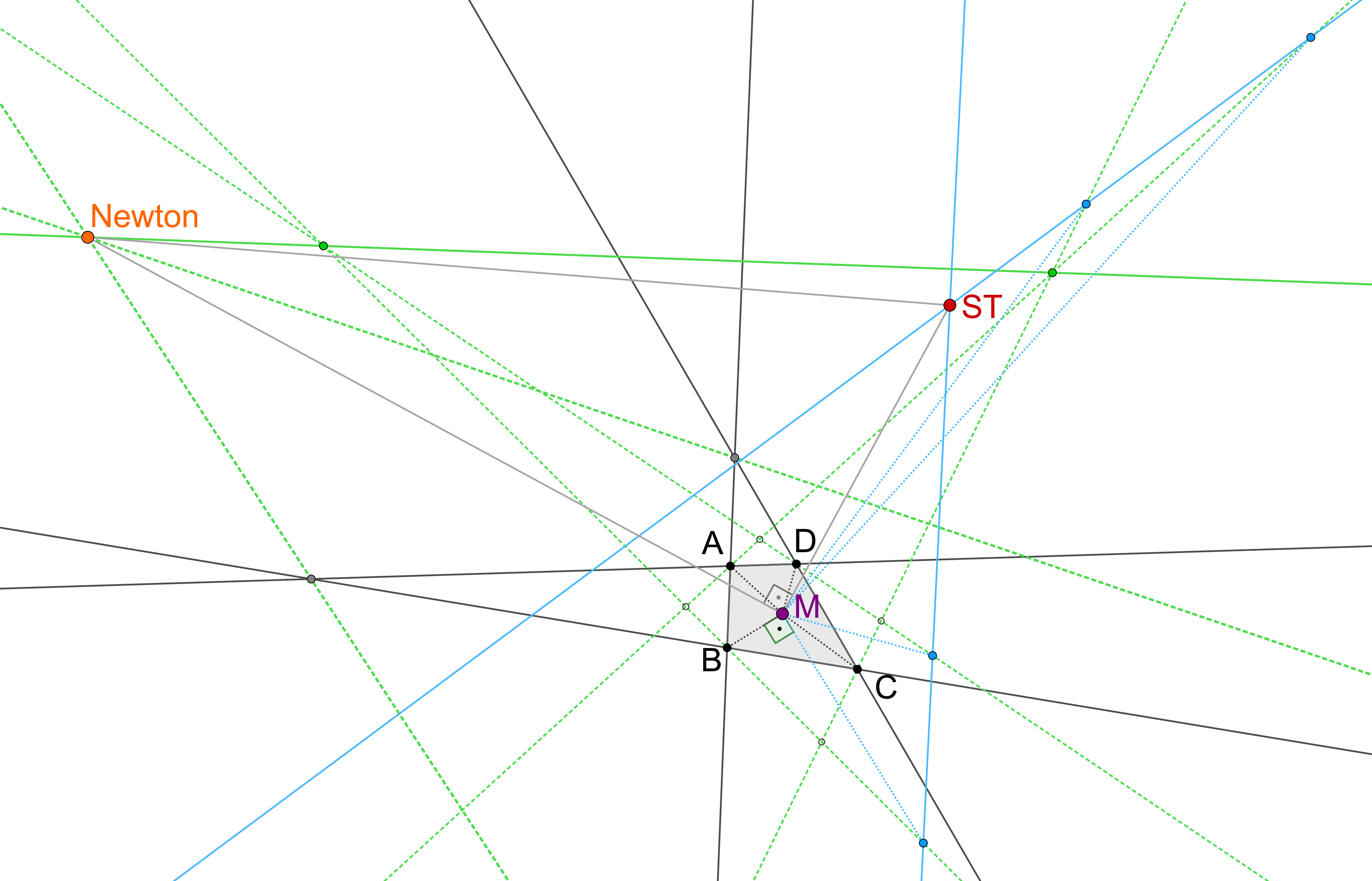}
\caption{Shatunov-Tokarev point and Newton point of a tetragon in a polar-Euclidean plane. $M$ is the absolute pole in this plane.}
\vspace*{0 mm}\end{figure}

\noindent\textit{Theorem} (L. Shatunov, A. Tokarev, cf. \cite{Sh}). Let $\mathbf{T} = ABCD$ be a tetragon  with vertices $A, B, C, D$ in a metric affine plane. As before, we imagine that this affine plane is embedded in a metric projective plane.\newpage
\noindent Then there exist a point $R$  with $[R,A]\cong [R,B]$ and $[R,C]\cong [R,D]$ and a point $S$ with $[S,A]\cong [S,D]$ and $[S,B]\cong [S,C]$.\\
If $\mathbf{T}$ is not cyclic, then the two points $R$ and $S$ are uniquely determined and distinct.  If $\mathbf{T}$ is not a parallelogram, it has a Newton line and the line $R{\,\vee\,}S$ is perpendicular to this Newton line. The line $R{\,\vee\,}S$ is now called the \textit{Shatunov-Tokarev line}. \vspace*{0.0 mm}\\
If $\mathbf{T}$ is a parallelogram but not cyclic, the Shatunov-Tokarev line is the line at infinity.\\
If $\mathbf{T}$ is a cyclic tetragon, then the center $M$ of the circumcircle of $\mathbf{T}$ is a point with $[M,A]\cong [M,B] \cong [M,C]\cong [M,D]$. In case $\mathbf{T}$ is not a parallelogram, we define the Shatunov-Tokarev line as the line through $M$ perpendicular to the Newton line.\vspace*{1.0 mm}\\
If $\mathbf{T}$ is not a parallelogram, the Shatunov-Tokarev line is the radical line of two special circles. One circle has its center at the midpoint $M_{\!AC}$ of $[A,C]$, the center of the other circle is the midpoint $M_{\!BD}$ of $[B,C]$. The radius of the first circle is the distance from $M_{\!BD}$ to $B$, the radius of the second the distance from $M_{\!AC}$ to $A$.\vspace*{-1.5 mm}\\

\begin{figure}[!htbp]
\begin{centering}
\includegraphics[height=5.9cm]{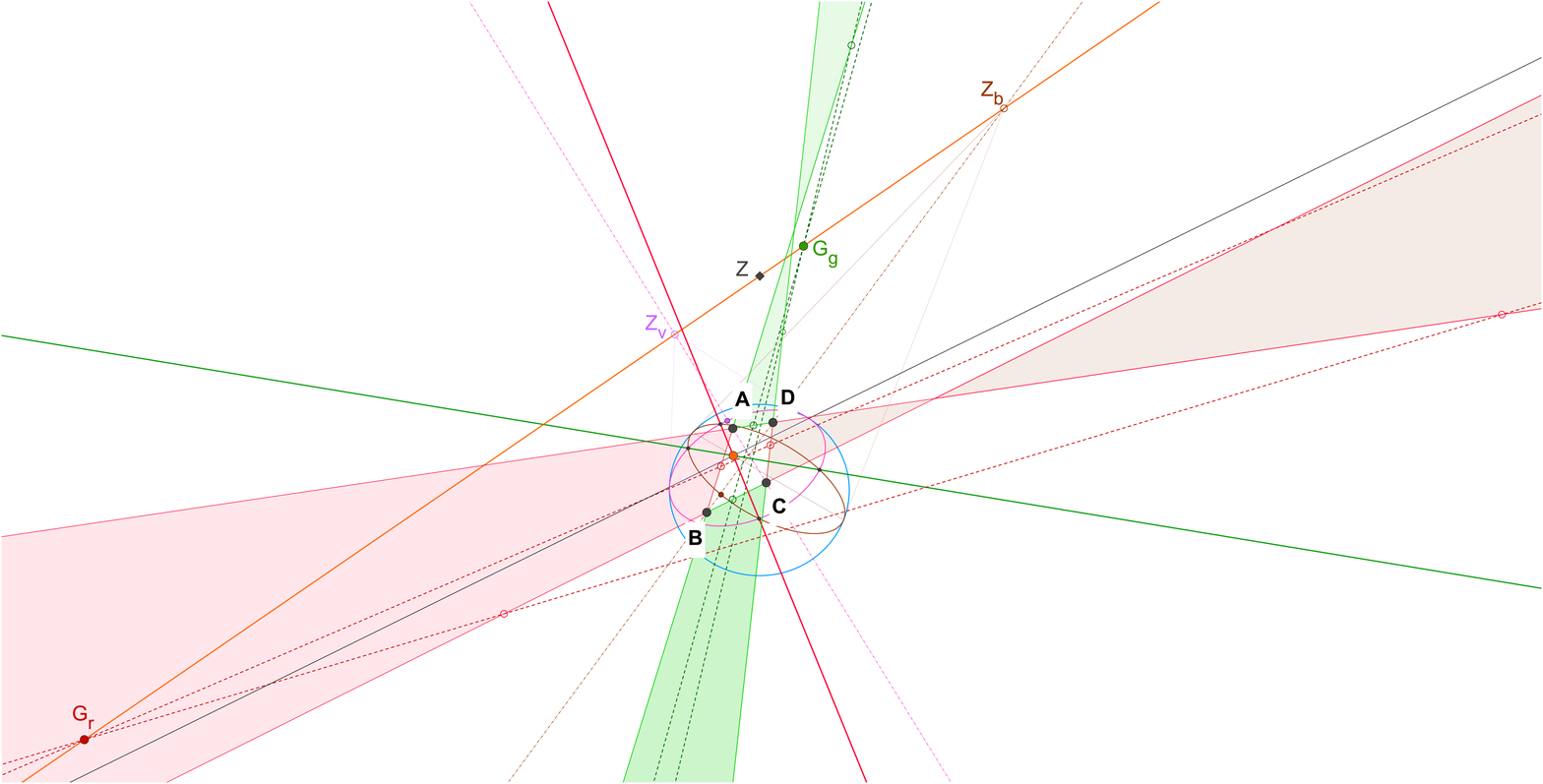}
\caption{Shatunov-Tokarev lines (red and green) and Newton line (orange) of a pair $(\mathbf{T},\mathbf{T^c})$ of tetragons with vertices $A,B,C,D$. The red tetragon is $\mathbf{T}^{(+-+-)}$. Both tetragons have proper centroids, $G_{\!r}$ and $G_{\!g}$. Besides these two points, the proper midpoint $Z_v$ of $[A,C]_-$, the proper midpoint $Z_b$ of $[B,D]_-$, and the proper midpoint $Z$ of $[Z_v,Z_b]_+$ are points on the Newton line. $G_{\!r}$, $Z_v$, $G_{\!g}$,$Z_b$ form a harmonic range.\\
Two circles are shown, one purple and one brown.
The purple circle is constructed as follows. Its center is $Z_v$, and the mirror images of points $B$ and $D$ in $Z$ are anisotropic points on the circle. The center of the brown circle is $Z_b$ and we can get an anisotropc point on the circle when we reflect point $A$ or point $C$ in $Z$. Both Shatunov-Tokarev lines are radical lines of the two circles.\\
\,}
\includegraphics[height=5.8cm]{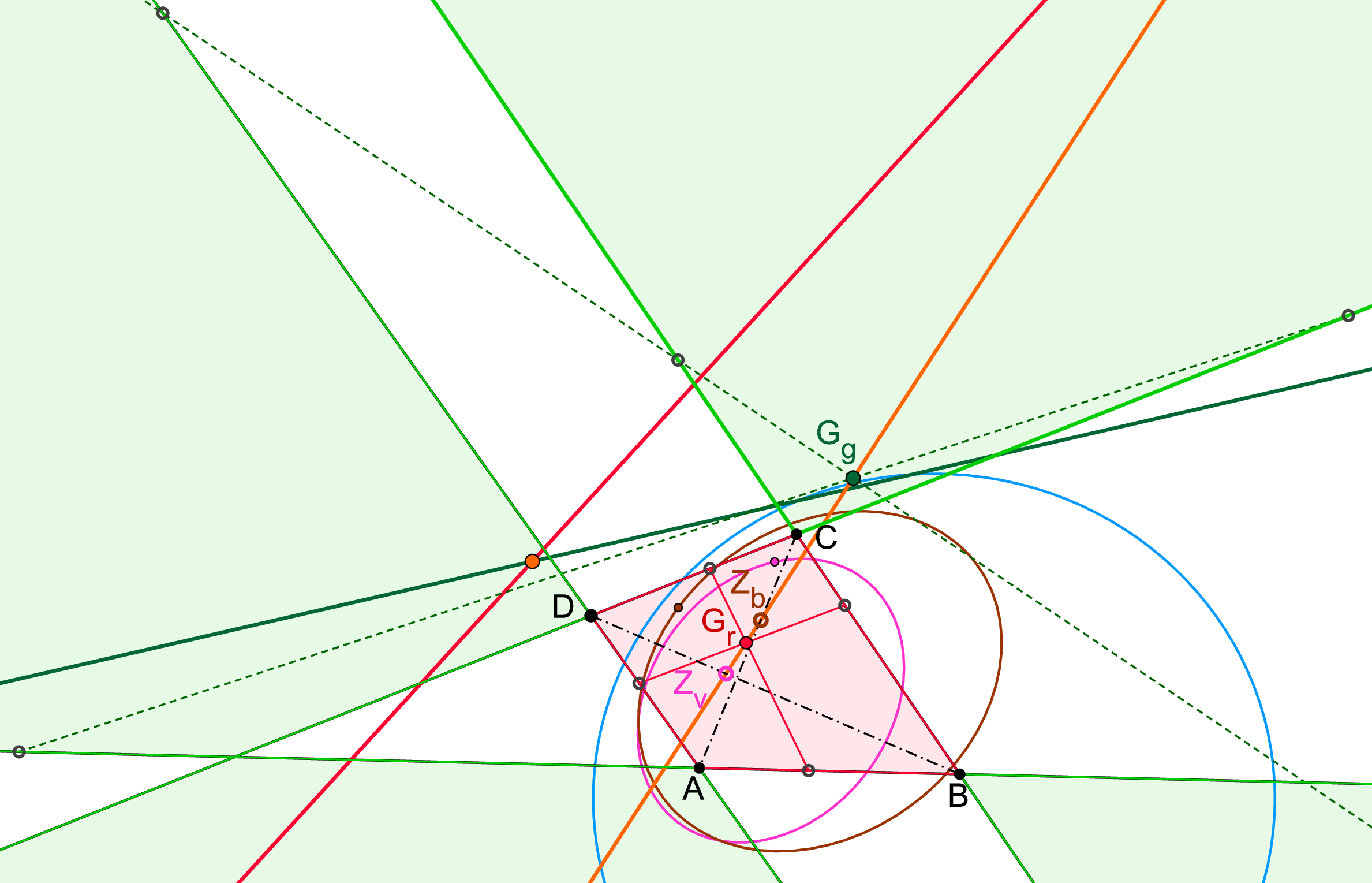}
\caption{Shatunov-Tokarev lines (red and dark green) and Newton line (orange) of a pair $(\mathbf{T},\mathbf{T^c})$ of tetragons in a hyperbolic plane. The red tetragon is $\mathbf{T}^{(++++)}$, the green $\mathbf{T}^{(----)}$. Both tetragons have pseudo-centroids. The purple circle and the brown circle are constructed as the ones in the previous figure. But the Shatunov-Tokarev lines are not radical lines of these circles. }
\end{centering}
\vspace*{-3.2 mm}\end{figure}

With some modifications, the Shatunov-Tokarev theorem can be applied to quadrilaterals in elliptic and hyperbolic planes. First, we introduce some new terms and notations. \vspace*{-1.5 mm} \\

A tetragon $\mathbf{T} = ABCD$ in a projective plane is not uniquely defined by the specification and order of its vertices. There are actually sixteen different such tetragons.
We can write each of these tetragons as the union of its sides, for example $[A,B]_+\cup [B,C]_-\cup [C,D]_- \cup [D,A]_+$. But we prefer to write this sequence more briefly $\mathbf{T^{(+--+)}}.$\\
All these tetragons share the same triple  $P_1,P_2,P_3$ of diagonal points.\\
For each tetragon $\mathbf{T}$ there is exactly one other tetragon $\mathbf{T}^c$ which has the same vertices but completely different sides; this we call the complementary of $\mathbf{T}$. For example, $\mathbf{T}=\mathbf{T^{(+--+)}}$ and $\mathbf{T}^c=\mathbf{T^{(-++-)}}$ are complementary to each other.\vspace*{-1.5 mm}\\

In an elliptic and in a hyperbolic plane, each of the tetragons with vertices $A,B,C,D$  has its own \textit{semi-centroid}; it is the intersection of the two lines which connect the inner semi-midpoints of opposite sides of the tetragon.\\
Definition: The \textit{Shatunov-Tokarev line} of $\mathbf{T}$ is the dual of the semi-centroid of $\mathbf{T}^c$. (Thus, the Shatunov-Tokarev line of $\mathbf{T}^c$ is the dual of the semi-centroid of $\mathbf{T}$.) \\
\textit{Remark}: If $\{A\}\cong\{B\}\cong\{C\}\cong \{D\}$, then there is at least one point $R$ on this Shatunov-Tokarev line with $[A,R]_+\cong [B,R]_+$ and $[C,R]_+\cong [D,R]_+$ and at least one point $S$ on this line  with $[A,S]_+\cong [D,S]_+$ and $[B,S]_+\cong [C,R]_+$. \vspace*{0.8 mm}\\
The Shatunov-Tokarev lines of $\mathbf{T}$ and $\mathbf{T^c}$ meet at some point $P$. We call the dual line of this point \textit{Newton line} of $\mathbf{T}$ and $\mathbf{T^c}$. Thus, the \textit{Newton line} of $\mathbf{T}$ is the line through the semi-centroid of $\mathbf{T}$ and the semi-centroid of $\mathbf{T}^c$.\vspace*{-1 mm} \\

We choose this definition of Newton line because: (1) Not for every tetragon the semi-centroid and the two semi-midpoints of the diagonals lying within this tetrahedron are collinear, therefore the definition of the Newton line of a tetrahedron in a metric affine plane cannot be transferred 1:1 to a tetrahedron in an elliptic or a hyperbolic plane.
(2) If there is an even number of plus signs in the exponent of $\mathbf{T}^{(\cdots\cdot)}$, then one semi-midpoint of $[A,C]_+$ (the inner or the outer) and one semi-midpoint of $[B,D]_+$ lie on this Newton line. \textit{Proof} of (2). We take $P_1=[1,0,0], P_2=[0,1,0], P_3=[0,0,1]$ as a generating system for the plane. There are nonzero real numbers $s,t$  and positive real numbers $a,b,c,d$ such that 
$ A^\circ = ({-}s,t,1)/a{\,,\,}B^\circ =({-}s,{-}t,1)/{b}{\,,\,}$ $C^\circ = (s,{-}t,1)/{c}{\,,\,}D^\circ$ $ = (s,t,1)/{d}\,$. \\
Further calculations depend on which quadrugon $\mathbf{T}$ represents. We do the calculations for $\mathbf{T}=\mathbf{T}^{(+--+)}$; calculations for the remaining cases are quite similar.\\
In this case, $\{\,[p_1{:}p_2{:}p_3]\,|\, p_1 t (a b{\,+\,}c d)+p_2 s (a d{\,+\,}b c)+p_3s\,t (a{\,+\,}c) (d{\,-\,}b)=0\,\}$ is the Newton line of $\mathbf{T}$ and of $\mathbf{T^c}$. Now it can be easily checked that $A^\circ\!+C^\circ$ and $B^\circ\!-D^\circ$ lie on the Newton line while $A^\circ\!-C^\circ$ and $B^\circ\!+D^\circ$ do not.\;$\Box$\vspace*{1.7 mm} \

Without giving a proof, we state that neither a semi-midpoint of $[A,C]_+$ nor a semi-midpoint of $[B,D]_+$ can lie on the Newton line if the number of plus signs is uneven.\\

\noindent \textit{Theorem} 6. If $\mathbf{T}=\mathbf{T}^{(.\dots)}$ is a tetragon with vertices $A,B,C,D$ in an elliptic or a hyperbolic plane. We call the line which passes through the semi-centroids of $\mathbf{T}$ and $\mathbf{T^c}$ the \textit{Newton line} of $\mathbf{T}$. If the number of plus signs in the exponent of $\mathbf{T}$ is even, then the Newton line passes through one semi-midpoint (inner or outer) of $[A,C]_+$ and one semi-midpoint of $[B,D]_+$.

%Two distinct tetragons with vertices $A,B,C,D$ have distinct semi-centroids.\\
%One of the two semi-midpoints of $[A,C]_+$ is a point inside the tetragon, the other lies outside, and the same applies to $[B,D]_+$.\vspace*{-1.5 mm}\\

\section*{Anne's theorem}

\noindent \textit{Theorem}  (Anne's Theorem for the Euclidean plane)\footnote{$^)$Pierre-Leon Anne (1806–1850)}$^)$ .
Let $\mathbf{T}$ be a convex tetragon in the Euclidean plane with vertices $A,B,C,D$. We assume  
that $\mathbf{T}$ is not a parallelogram, so it has a Newton line. Let $Q$ be a point inside the tetragon.\\
Then the sum of the areas of the triangles $ABQ$ and $CDQ$ is equal to the sum of the areas of the triangles $BCQ$ and $DAQ$ if and only if $Q$ lies on the Newton line. \vspace*{0.8 mm}\\
For a proof of this theorem see for example \cite{HS}. This paper also provides some information about the history of Anne's and Newton's theorems.\vspace*{-1 mm}\\

%\begin{figure}[!htbp]
%\includegraphics[height=5cm]{Anne.eps}
%\caption{}
%\vspace*{1 mm}\end{figure}

This theorem does not apply in elliptic and hyperbolic planes if we take the angle excess as a measure for the area of a triangle, as can easily be seen by the following example.\vspace*{-1.6 mm}\\

Let $A, B, C$ be vertices of a triangle with sides $[A,B]_+, [B,C]_+, [C,A]_+$.
Let us assume that the segments $[A,C]_+$ and $[B,C]_+$ are not congruent. In this case, the line through $B$ that divides the triangle into two parts of equal area, intersects the side $[C,A]_+$ at a point which differs from the midpoint $M$ of this side, see \cite{Ev1}. The subtriangles $AMB$ and $BMC$ have different areas.\\
Let $D$ be the image of $B$ under a reflection in $M$. Under this reflection, the triangles $ABC$, $CDA$, $CDA$ are mapped onto the equal-area triangles $CDA$, $CMD$, $MAD$, respectively.
Thus, the sum of the areas of the triangles $ABF$ and $CDF$ is different from the sum of the areas of the triangles $BCF$ and $DAF$. However, the point $M$ lies on the Newton line of the tetragon $\mathbf{T}^{(++++)}$ with vertices $A,B,C,D$. $\;\;\;\;\Box$\vspace*{0.8 mm}\\
\textit{Remark}. We know, of course, that in Euclidean geometry the angle excess is not a suitable measure for the area.\vspace*{0.8 mm}\\

There is a function in elliptic and hyperbolic geometry that shares some characteristic properties with the area function. It is called \textit{Staudtian} and is defined as follows.\vspace*{0.5 mm}\\
We can find points $P_1, P_2, P_3$ in both the elliptic and hyperbolic plane, such that the metric tensor {$\mathfrak{G}$ has a diagonal form with respect to these points, and we can assume that in the elliptic case $\mathfrak{g}_{11}=\mathfrak{g}_{22}=\mathfrak{g}_{33}=1$, in the hyperbolic case $\mathfrak{g}_{11}=\mathfrak{g}_{22}=1, \mathfrak{g}_{33}=-1$.\\
We define a function $\psi$ on the powerset of $\textrm{P}(\mathbb{R}^3)$ with values in $\mathbb{R}$.  If $\mathcal{S}$ is a subset of $\textrm{P}(\mathbb{R}^3)$, then $\psi(\mathcal{S})$ = 0, unless $\mathcal{S}$ consists of a finite number of anisotropic points. If $\mathcal{S}$ consists of $s$ 
anisotropic points $Q_1,\dots, Q_s$, then  $\psi(\mathcal{S}) = \det((Q_i^\circ{\,\scriptscriptstyle{[\mathfrak{G}]}\,}Q_j^\circ)_{1\le i,j\le s})$.\\
If $\psi(\mathcal{S})\ne 0$, we call the number $\sigma(\mathcal{S}){\,{{\scriptstyle{:}}\!\!=}\,}\frac{1}{(s-1)!}\sqrt{|\psi(\mathcal{S})|}$ the \textit{Staudtian} of $\mathcal{S}$, cf. \cite{Ho,Ev1,Ev2}.\\
The function $\sigma$ has following properties:\\
(1) Given nine points $A_1,A_2,A_3,B_1,B_2,B_3,C_1,C_2,C_3$ such that $\sigma(\{A_1,B_1,C_1\})=\;\;$ $\sigma(\{A_2,B_2,C_2\})$ and $\sigma(\{A_2,B_2,C_2\})=$ $\sigma(\{A_3,B_3,C_3\})$.\\Then $\sigma(\{A_1,B_1,C_1\})=\,$ $\sigma(\{A_3,B_3,C_3\})$.\\
(2) Let $A, B, C$ be three noncollinear anisotropic points. Let $P$ be the dual of the line $A\vee B$, and let us assume that $P$ is an anisotropic point distinct from $C$.\\
(2a) Then the line $C\vee P$ intersects the line $A{\,\vee\,}B$ in a point we call $Q$, and $\sigma(\{A,B,C\}) = \frac{1}{2}\,\sigma(\{A,B\})\,\sigma(\{C,Q\}).$\\  
(2b) Let $\tilde{C}$ be yet another anisotropic point. Then $\sigma(\{A,B,C\}) = \sigma(\{A,B,\tilde{C}\})$ iff $\tilde{C}$ is a point on the circle around $P$ through $C$.\\
(2c) If $Q$ is any anisotropic point on $A{\,\vee\,}B$, then  $\sigma(\{A,Q,C\})=\sigma(\{Q,B,C\})$ iff $Q$ is an inner or outer midpoint of $[A,B]_+$.\\ 
(2d) If $D$ be a point on $A\vee B$, then\\
\centerline{$\sigma(\{A,D,C\}) + \sigma(\{D,B,C\}) = \sigma(\{A,B,C\})\;\;\;\;\;\mathrm{if}\!\mathrm{f}\;\;\;\;\;D=A\,\vee\,D=B$.}\\ 
In contrast to the area, the Staudtian lacks the property of additivity.\\

We present three formulas on which we rely in the proof of the next theorem.\\
If $\boldsymbol{v}_1,\boldsymbol{v}_2,\boldsymbol{v}_3,\boldsymbol{v}_4,\boldsymbol{w} \in \mathbb{R}^3$, then\\
(1) $\det((\boldsymbol{v}_{i}{\,\scriptscriptstyle{[\mathfrak{G}]}\,}\boldsymbol{v}_{j})_{1\le i,j\le 3}) = (\det(\boldsymbol{v}_1,\boldsymbol{v}_2,\boldsymbol{v}_3))^2$\\
(2) $\det(\boldsymbol{w},\boldsymbol{v}_1,\boldsymbol{v}_2)-\det(\boldsymbol{w},\boldsymbol{v}_2,\boldsymbol{v}_3)+\det(\boldsymbol{w},\boldsymbol{v}_3,\boldsymbol{v}_4)-\det(\boldsymbol{w},\boldsymbol{v}_4,\boldsymbol{v}_1) = \det(\boldsymbol{w},\boldsymbol{v}_1+\boldsymbol{v}_3,\boldsymbol{v}_2+\boldsymbol{v}_4)$\\
(3) Given a triangle $ABC$ in an elliptic or in a hyperbolic plane, then
$\sigma(\{A,B,C\})= \frac{1}{2} s \det(A^\circ,B^\circ,D^\circ)$, where $s=1$ or $s=-1$, depending on whether the triangle is positively or negatively oriented.\\

%$|\det(\boldsymbol{w},\boldsymbol{v}_1,\boldsymbol{v}_2)|-|\det(\boldsymbol{w},\boldsymbol{v}_2,\boldsymbol{v}_3)|+|\det(\boldsymbol{w},\boldsymbol{v}_3,\boldsymbol{v}_4)|-|\det(\boldsymbol{w},\boldsymbol{v}_4,\boldsymbol{v}_1))| = |\det(\boldsymbol{w},\boldsymbol{v}_1+\boldsymbol{v}_3,\boldsymbol{v}_2+\boldsymbol{v}_4)|$

\noindent\textit{Theorem} 7. Let $\mathbf{T}=ABCD$ be a strictly convex tetragon with vertices $A, B, C, D$ in an elliptic or in a hyperbolic plane. Let us assume that the boundary of $\mathbf{T}$ is formed by the segments $[A,B]_+,[B,C]_+,[C,D]_+$ and $[D,A]_+$ and let us further assume that in the hyperbolic case the vertices are all inside or all outside the absolute conic.\\
If $P$ is an anisotropic point inside the tetragon, then $\sigma(P,A,B)+\sigma(P,C,D)=\sigma(P,B,C)+\sigma(P,D,A)$
if and only if $P$ is a point on the Newton line.\vspace*{0.6 mm}\\
\textit{Proof}. Let $M_{AC}$ and $M_{BD}$ be the inner midpoints of $[A,C]_+$ and $[B,D]_+$, respectively. Put $s{\,{{\scriptstyle{:}}\!\!=}\,}\mathrm{sgn}(\det(P^\circ,A^\circ,B^\circ))$. The triangles $PAB, PBC,PCD,PDA$ all have the same orientation, therefore\\
$\hspace*{10 mm}0\; = 2 s\,\big(\sigma(\{P,A,B\})-\sigma(\{P,B,C\})+\sigma(\{P,C,D\})-\sigma(\{P,D,A\})\big)\vspace*{0.4 mm} \\
\hspace*{12.6 mm} = \det(P^\circ,A^\circ,B^\circ)-\det(P^\circ,B^\circ,C^\circ)+\det(P^\circ,C^\circ,D^\circ)-\det(P^\circ,D^\circ,A^\circ)\vspace*{0.7 mm}$\\
$\displaystyle\hspace*{13.5 mm} = \det(P^\circ,{A^\circ+C^\circ},{B^\circ+D^\circ})\;.\;\; {\hspace*{56 mm}} \Box$ \\

% {\,\smallsetminus}{\hspace{-4.3pt}\smallsetminus\,}
%{\,\smallsetminus}{\hspace{-4.3pt}\smallsetminus\,}
%\noindent\textit{Theorem} 11? (Euler's quadrilateral law for the elliptic and the hyperbolic plane)\\ 
%The same assumptions apply as in the previous theorem. Then,\\
%{\hspace*{11 mm}$\sigma^2(\{A,B\})+\sigma^2(\{B,C\})+\sigma^2(\{C,D\})+\sigma^2(\{D,A\})$\\
%${\hspace*{6 mm} = {\hspace*{1.6 mm} \sigma^2(\{A,C\})+\sigma^2(\{B,D\})+4\,\sigma^2(\{M_{A,C},M_{B,D}\})}}$. 

\end{document}